\documentclass[review]{elsarticle}

\usepackage{lineno,hyperref}
\modulolinenumbers[5]

\journal{}
\usepackage[english,activeacute]{babel}
\usepackage{amsfonts}
\usepackage{amsmath}
\usepackage{url}
\usepackage{hyperref}
\usepackage{algorithm}
\usepackage{algpseudocode}
\usepackage{appendix}
\usepackage{empheq}
\usepackage{booktabs}
\usepackage{multirow}
\usepackage{mathtools}
\usepackage{tikz,pgfplots,pgfplotstable,booktabs}
\usepackage{xcolor}
\usepackage{bm}
\usepackage{enumitem}
\usepackage{lscape}
\usepackage{rotating}
\usepackage{amssymb}
\usepackage{subcaption}
\usepackage{graphicx}
\usepackage{floatpag}
\usepackage{xr}
\usepackage[maxfloats=40]{morefloats}
\newcommand{\myreferences}{Bibliography_all}

\usetikzlibrary{shapes,arrows,fit,backgrounds}
\tikzstyle{block} = [draw, rectangle, text width=3cm, text centered, minimum height=1.2cm, node distance=3.5cm,fill=white]
\tikzstyle{container} = [draw, rectangle, inner sep=0.2cm, fill=gray!50,minimum height=3.2cm]
\def\bottom#1#2{\hbox{\vbox to #1{\vfill\hbox{#2}}}}
\tikzset{
  mybackground/.style={execute at end picture={
      \begin{scope}[on background layer]
        \node[] at (current bounding box.north){\bottom{1cm} #1};
        \end{scope}
    }},
}
\usepackage[draft]{todonotes} 
\usetikzlibrary{calc}
\graphicspath{{Figures/}}

\begin{document}
%\begin{center}
%\Large
%    Essential Title Page Information
%\end{center}
%\begin{itemize}
%    \item \textbf{Title of the paper:}  Warm-starting constraint generation for mixed-integer optimization: A Machine Learning approach.
%    \item \textbf{Full name and email address of every contributing author:}
%    \begin{itemize}
%        \item Asunci\'on Jim\'enez-Cordero (asuncionjc@uma.es). Corresponding\newline author.
%        \item Juan Miguel Morales (juanmi82mg@gmail.com).
%        \item Salvador Pineda (spinedamorente@gmail.com).
%    \end{itemize}
%    \item \textbf{Full affiliation with country of all the contributing authors:} OASYS research group. Ada Byron research building. Street ``Arquitecto Francisco Pe\~nalosa'', 18, 29010. University of M\'alaga, M\'alaga, Spain.
%\end{itemize}
%\newpage
\begin{frontmatter}
    
%\title{Machine-learning-aided warm-start of\\ constraint generation methods for\\ online mixed-integer optimization}
\title{\textcolor{black}{Warm-starting constraint generation \\ for mixed-integer optimization:\\ A Machine Learning approach}}

%% Group authors per affiliation:
\author{Asunci\'on Jim\'enez-Cordero\corref{mycorrespondingauthor}}
\cortext[mycorrespondingauthor]{Corresponding author}
\ead{asuncionjc@uma.es}
\author{Juan Miguel Morales}
\ead{juanmi82mg@gmail.com}
\author{Salvador Pineda}
\ead{spinedamorente@gmail.com}

\address{OASYS Group. Ada Byron Research Building,\\ Arquitecto Francisco Pe\~nalosa St., 18, 29010,\\ University of M\'alaga, M\'alaga, Spain\\ Phone: + 34 951 95 29 25}

\begin{abstract}
Mixed Integer Linear Programs (MILP) are well known to be NP-hard (Non-deterministic Polynomial-time hard) problems in general.
Even though pure optimization-based methods, such as constraint generation, are guaranteed to provide an optimal solution if enough time is given, their use in online applications remains a great challenge due to their usual excessive time requirements. To alleviate their computational burden, some machine learning techniques (ML) have been proposed in the literature, using the information provided by previously solved MILP instances. Unfortunately, these techniques report a non-negligible percentage of infeasible or suboptimal instances.

By linking mathematical optimization and machine learning, this paper proposes a novel approach that speeds up the traditional constraint generation method, preserving feasibility and optimality guarantees. In particular, we first identify offline the so-called invariant constraint set of past MILP instances. We then train (also offline) a machine learning method to learn an invariant constraint set as a function of the problem parameters of each instance. Next, we predict online an invariant constraint set of the new unseen MILP application and use it to initialize the constraint generation method. This warm-started strategy significantly reduces the number of iterations to reach optimality, and therefore, the computational burden to solve online each MILP problem is significantly reduced. Very importantly, \textcolor{black}{all the feasibility and optimality theoretical guarantees of the traditional constraint generation method are inherited by our proposed methodology}. The computational performance of the proposed approach is quantified through synthetic and real-life MILP applications.

\end{abstract}


\end{frontmatter}

%\linenumbers

\section{Introduction}\label{sec: Introduction}

\textcolor{black}{Recent papers, e.g., \cite{gambella2021optimization, lin2017multi, yang2020on, zamfirache2022policy}, have shown the potential in combining  Mathematical Optimization and Machine Learning. Particularly,} Mixed Integer Linear Programming (MILP) is \textcolor{black}{known to be} a powerful and flexible tool for modeling and solving a wide variety of decision-making problems, as can be confirmed from \cite{conforti2014integer, wolsey2008mixed}. However, most MILPs are known to be \emph{NP-hard} (Non-deterministic Polynomial-time hard) and therefore, using them for online applications becomes challenging. For this reason, the design of novel Machine-Learning-assisted (ML) techniques that reduce the computational burden of MILPs has recently become a popular research topic. The works in \cite{bengio2021machine, bertsimas2021voice} are a valid proof of that. Although different learning strategies can be considered, all methods assume that a set of slight variations of the same problem have been previously solved, and their input data and solutions are available. This is a reasonable assumption since many optimization problems are frequently solved for a range of input parameters in online applications.

Several works have proposed approaches that seek to preserve optimality guarantees while substantially reducing the solution time. For instance, the authors of \cite{chmiela2021learning} design a data-driven methodology to improve the use of heuristics in branch-and-bound. To be more precise, the authors train a machine learning model that provides a schedule of heuristics, specifying when and for how long each heuristic is executed. They numerically show that this smart schedule is very likely to substantially diminish the time to solve the mixed-integer program to optimality. Even though the proposed method yields impressive results, collecting data from the branch-and-bound heuristics may be a difficult and solver-dependent task. Hence, alternative approaches that only need and learn from the information provided by the MILP optimal solution are also considered. One of the most commonly used strategies in this line consists in building, from this information, a \emph{simpler} formulation of the original MILP that is faster to solve. 

A pure optimization-based method that can be applied to this end is \emph{constraint generation}. A detailed explanation of this methodology can be found in \cite{minoux1989networks}. Essentially, it sequentially adds violated constraints until the optimal solution is found, at the expense of solving a possibly large number of MILPs. As a consequence, the computational burden associated with this strategy may be unacceptable for online applications. One way of alleviating such computational effort would be to provide a \emph{good} warm-start. That is, if a good initial set of constraints is given, then just a few iterations of the constraint generation method are executed, reducing the computational time. In this vein, the authors of \cite{pineda2020data} propose a learning strategy that uses a modified nearest neighbor methodology to screen out superfluous constraints of the Unit Commitment (UC) problem. The reader is referred to \cite{gomez2018electric} for further details about this problem. Similarly, the method presented in \cite{xavier2020learning} selects a subset of constraints for the same problem. All these works execute the constraint generation algorithm with an initial set of constraints that is inferred only from the constraints that appeared binding in previous instances of the problem, that is, from constraints that held with equality at the optimal solution. Using these binding constraints is not enough to recover the optimal solution in MILPs, as shown in \cite{pineda2020data}. Indeed, when integer variables are involved in the optimization problem, the optimal solution must satisfy the constraints included in the so-called \emph{invariant constraint set} defined in \cite{calafiore2010random}. Apart from the constraints that hold with equality, such an invariant constraint set includes other non-binding constraints that are crucial to attaining the optimal solution.

\textcolor{black}{Therefore, the contributions and objectives achieved in this paper are:
\begin{itemize}
    \item[-] The authors in \cite{pineda2020data} completely ignore the critical non-binding constraints. Hence, the learned warm-start set is not good enough and, to attain optimality, the constraint generation method may require a possibly large number of iterations, thus resulting in quite modest computational savings. In contrast, we show that effectively warm-starting the constraint generation method is essential to reaching the optimal solution of a MILP after just a few iterations of the algorithm.
    \item[-] To this aim, we identify offline an invariant constraint set for each past instance. Then, we train a machine learning model of our choice that returns the prediction of an invariant constraint set from MILP parameters.
    \item[-] Such a learned set is used to warm-start the constraint generation method of a new MILP. This way, the optimal solution is attained after running just a small number of iterations in the constraint generation method, and the computational time of the online MILP is considerably reduced, as empirically shown using synthetic and real-life instances.''
\end{itemize}
}

% The data-driven methods proposed in \cite{pineda2020data} and \cite{xavier2020learning} completely ignore such non-binding constraints. Therefore, the learned warm-start set is not good enough, and, to attain optimality, the constraint generation method may require a possibly large number of iterations, thus resulting in quite modest computational savings.

%In this paper, we show that effectively warm-starting the constraint generation method is essential to reaching the optimal solution of an MILP in a short period of time. To this aim, we identify offline an invariant constraint set for each past instance. Then, we train a machine learning model of our choice that returns the prediction of an invariant constraint set from MILP parameters. Such a learned set is used to warm-start the constraint generation method of a new MILP. This way, the optimal solution is attained after running just a small number of iterations in the constraint generation method, and the computational time of the online MILP is considerably reduced, as empirically shown using synthetic and real-life instances. \textcolor{red}{[Asun: Comment R3.3.]}

The remainder of this paper is structured as follows: Section \ref{sec: Invariant Constraint Set in MILPs} presents MILP notation and discusses the difficulties of computing an invariant constraint set when integer variables appear. Section \ref{sec: Methodology} details the proposed methodology. Section \ref{sec: Computational Experience} focuses on the numerical experiments. Finally, Section \ref{sec: Conclusions} provides the conclusions and future research lines.

\section{Invariant Constraint Set in MILPs}\label{sec: Invariant Constraint Set in MILPs}

A general MILP can be formulated as follows:

    \begin{subequations}%\label{subeq: MILP}
 \makeatletter
        \def\@currentlabel{$P_{\bm{\theta}}[\mathcal{J}]$}
        \makeatother
        \label{eq:w}
         \renewcommand{\theequation}{$P_{\bm{\theta}}[\mathcal{J}]$}
  \begin{empheq}[left=\empheqlbrace]{align}
    \min\limits_{\boldsymbol{z}\in\mathbb{R}^{n} \times \mathbb{Z}^q}&\boldsymbol{c}^{\intercal}\boldsymbol{z}\label{subeq: MILP}  \\
    \text{s.t. }  & \boldsymbol{a}_j^{\intercal}\boldsymbol{z} \leq b_j, \quad \forall j\in \mathcal{J}\notag 
    \end{empheq}
  \end{subequations}
\noindent where $\boldsymbol{c},\boldsymbol{a}_j\in\mathbb{R}^{n+q}, \forall j$, $b_j\in\mathbb{R},\, \forall j$, are input parameters  and $\boldsymbol{z}= (\boldsymbol{x}, \boldsymbol{y})$ is the decision variable vector formed by the continuous variables $\boldsymbol{x}\in\mathbb{R}^n$ and the integer variables $\boldsymbol{y}\in\mathbb{Z}^q$.  For convenience, we collect all those input parameters into the set $\bm{\theta}$, that is, $\bm{\theta} = \{\boldsymbol{c}, \boldsymbol{a}_j, b_j, \forall j \in \mathcal{J} \}$, so that \eqref{subeq: MILP} denotes the optimization problem with the set of input parameters $\bm{\theta}$ and the set of constraints $\mathcal{J}$. Note that there is some abuse of notation in \eqref{subeq: MILP}. As it is written, $\mathcal{J}$ and $\boldsymbol{\theta}$ can change independently. However, this is not true, since the definition of $\boldsymbol{\theta}$ explicitly depends on $\mathcal{J}$. Hence, to be rigorous, we should write $(P_{\boldsymbol{\theta}_{\mathcal{J}}}[\mathcal{J}])$. Nevertheless, in order to make the notation clearer, we remove the subindex $\mathcal{J}$ in $\boldsymbol{\theta}_{\mathcal{J}}$. In addition, for simplicity, we assume that problem \eqref{subeq: MILP} is bounded and feasible, and that its optimal solution $\boldsymbol{z}^*_{\bm{\theta}}[\mathcal{J}]$ is assumed to be unique. Note, however, that if multiple optimal solutions appear, retaining just one of them is enough for our proposal.

The feasible region defined by constraints in $\mathcal{J}$ includes the subset of the so-called \emph{binding constraints} $\mathcal{B}$. Particularly, $\mathcal{B}$ is comprised of the inequality constraints that hold with equality at the optimal solution, i.e., $\mathcal{B} = \{j\in\mathcal{J}:\boldsymbol{a}_j^{\intercal}\boldsymbol{z}^*_{\theta}[\mathcal{J}] = b_j\}$. Besides, according to \cite{calafiore2010random}, a subset of constraints $\mathcal{S}\subset \mathcal{J}$ is defined to be an \emph{invariant constraint set}, if the objective values of problems $(P_{\bm{\theta}}[\mathcal{J}])$ and $(P_{\bm{\theta}}[\mathcal{S}])$ coincide, i.e., if $\boldsymbol{c}^{\intercal}\boldsymbol{
z}^*_{\bm{\theta}}[\mathcal{J}]=\boldsymbol{c}^{\intercal}\boldsymbol{
z}^*_{\bm{\theta}}[\mathcal{S}]$.
Note that, following the previous definition, a unique invariant constraint set may not exist. The relationship between these two sets of constraints, $\mathcal{B}$ and $\mathcal{S}$, depends on whether problem \eqref{subeq: MILP} includes integer variables or not. Therefore, we discuss first the case in which all variables are continuous and then the more general case with both continuous and integer variables.

Let us first assume a particular case of \eqref{subeq: MILP} that only includes continuous variables, i.e., $q = 0$ and $\boldsymbol{z} = \boldsymbol{x}\in\mathbb{R}^n$. Using the optimal solution, $\boldsymbol{z}^*_{\bm{\theta}}[\mathcal{J}]$, it is straightforward to determine the set of binding constraints $\mathcal{B}$. Actually, the invariant constraint sets $\mathcal{S}$ of linear programming problems must contain all the constraints in $\mathcal{B}$, as the authors of \cite{bertsimas1997introduction} affirm. In other words, when the decision variables of the optimization problem are continuous, one can choose $\mathcal{S} = \mathcal{B}$. This way, instead of solving the original optimization problem \eqref{eq:w}, which may involve a high computational cost, the optimal solution is computed through the reduced problem $\left(P_{\bm{\theta}}[\mathcal{S}]\right) = \left(P_{\bm{\theta}}[\mathcal{B}]\right)$. Such a reduced problem includes fewer constraints, typically implies a lower computational effort and is thus more appropriate for online applications. 

Now take the more general case in which problem \eqref{subeq: MILP} includes both continuous and integer variables. In this case, authors of \cite{pineda2020data} demonstrate through an illustrative example that an invariant constraint set $\mathcal{S}$ can also include constraints that are non-binding at the optimum. These constraints play a critical role when solving MILPs since they cannot be removed from the original feasible region without impairing the feasibility, and thus the optimality, of the so-obtained solution.  In other words, when integer variables are involved in the optimization problem, an invariant constraint set includes not only binding constraints but also some of the non-binding ones, i.e., it holds that $\mathcal{S} \supset \mathcal{B}$. Like in the previous case, assume we have access to the optimal solution $\boldsymbol{z}^*_{\bm{\theta}}[\mathcal{J}]$. While the set of binding constraints can also be easily determined by evaluating all constraints at the optimum, identifying such a subset of critical non-binding constraints is a more challenging task when integer variables appear. In addition, ignoring these non-binding constraints is quite dangerous since the optimal solution of the reduced problem may violate some of the constraints in the set $\mathcal{J}$, increasing the number of iterations required by the constraint generation method.

Our claim is that, if tuples $(\boldsymbol{\theta}_t, \mathcal{S}_t)$ for previously solved instances $t = 1, \ldots, T$ are used to train any machine learning algorithm ML$(\cdot)$, then the solution of the reduced formulation of a new unseen instance $\tilde{t}$ with the predicted invariant constraint set, $P_{\theta_{\tilde{t}}}[\mathcal{S}_{\tilde{t}}]$ would more likely be feasible (and optimal) for the original problem $P_{\theta_{\tilde{t}}}[\mathcal{J}]$, than if only the tuples with the binding constraints $(\boldsymbol{\theta}_t, \mathcal{B}_t), \, \forall t$ are considered, as done in \cite{pineda2020data}. As a consequence, the number of iterations to be run online in the constraint generation method is decreased. Hence, it is desirable to warm-start the constraint generation algorithm with a constraint set as close to an invariant constraint set of the original MILP as possible. This way, the number of iterations executed by the algorithm is reduced, and so is its running time.

\section{Methodology}\label{sec: Methodology}

The goal of this paper is to develop a data-driven approach that guarantees the optimal solution of an MILP in a reduced computational time. To this aim, we propose a methodology that efficiently warm-start the constraint generation method so that the number of iterations executed (online) is as small as possible. The key point of our strategy is to initialize (warm-start) the constraint generation method using a predicted invariant constraint set, $\mathcal{S}\subset\mathcal{J}$, learned from past instances. 

 In particular, we propose a data-driven strategy that takes advantage of the information $(\boldsymbol{\theta}_t, \mathcal{S}_t)$ provided from the previously solved instances $t = 1, \ldots, T$ to predict an invariant constraint set, $\mathcal{S}_{\tilde{t}} = \text{ML}(\boldsymbol{\theta}_{\tilde{t}})$ for a new unseen problem $\tilde{t}$, using a machine learning model of our choice, $\text{ML}(\cdot)$. See \cite{friedman2001elements} for more details about the main machine learning tools. Then, once such a set $\mathcal{S}_{\tilde{t}}$ is predicted, we only need to run a smaller number of iterations than those required for the original constraint generation method to converge to the optimal solution of the MILP. Indeed, if the predicted invariant constraint set exactly coincides with the actual one, then, only one iteration of the constraint generation method needs to be executed. In other words, when the prediction is perfect, then it suffices to solve $P_{\boldsymbol{\theta}_{\tilde{t}}}[\mathcal{S}_{\tilde{t}}]$.
 
As mentioned in Section \ref{sec: Invariant Constraint Set in MILPs}, given an optimization problem, finding an invariant constraint set is a challenging task when integer variables appear. Developing an efficient procedure to learn an invariant constraint set in MILPs, $\mathcal{S}\subset \mathcal{J}$, is a relevant research question that has not yet been properly answered in the literature. We propose in this paper a methodology that aims at determining, for each training instance $t$, an invariant constraint set $\mathcal{S}_t$. The proposed approach to construct $\mathcal{S}_t$ for each instance $t$ is also based on a constraint generation procedure. To compute an invariant constraint set $\mathcal{S}_t$, for a previously solved instance $t$, we proceed as follows: we initialize the invariant constraint set $\mathcal{S}_t$ with the set of binding constraints $\mathcal{B}_t$. At each iteration, a reduced problem $P_{\bm{\theta}_t}[\mathcal{S}_t]$ is solved giving the optimal solution $\boldsymbol{z}^*_{\bm{\theta}_t}[\mathcal{S}_t]$. If all original constraints are satisfied, i.e., $\boldsymbol{a}_j^{\intercal}\boldsymbol{z}^*_{\bm{\theta}_t}[\mathcal{S}_t]\leq b_j,\, \forall j\in\mathcal{J}\setminus \mathcal{S}_t$, then the algorithm terminates. Otherwise, the most violated constraint is included in the set $\mathcal{S}_t$, and a new iteration is run. This algorithm is run offline for each training instance $t$. Hence, adding, at each iteration, only the most violated constraint is appropriate for our proposal since the online running time is not affected. However, alternative strategies, such as including in the set $\mathcal{S}_t$ all the violated constraints at each iteration, can be considered. The pseudocode of the proposed procedure is given in Algorithm \ref{alg: identifying support constraints}.

\begin{algorithm}[!htb]
  \begin{enumerate}[start=0,label={\arabic*)}]
    \item
     Initialize $\mathcal{S}_t = \mathcal{B}_t$.
    \item Solve $P_{\bm{\theta}_t}[\mathcal{S}_t]$ with solution $\boldsymbol{z}^*_{\bm{\theta}_t}[\mathcal{S}_t]$.
    \item If $\max\limits_{j\in\mathcal{J}\setminus\mathcal{S}_t} \left\{\boldsymbol{a}^{\intercal}_j \boldsymbol{z}^*_{\bm{\theta}_t}[\mathcal{S}_t] - b_j\right\} > 0$, go to step 3. Otherwise, stop.
    \item $\mathcal{S}_t := \mathcal{S}_t \cup \left\{ \arg \max \limits_{j\in\mathcal{J}\setminus\mathcal{S}_t} \{\boldsymbol{a}^{\intercal}_j \boldsymbol{z}^*_{\bm{\theta}_t}[\mathcal{S}_t] - b_j\}\right\}$, go to step 1.
  \end{enumerate}
   \caption{Identifying an invariant constraint set for each instance $t$}
   \label{alg: identifying support constraints}
\end{algorithm}

After running Algorithm \ref{alg: identifying support constraints}, the information $(\boldsymbol{\theta}_t, \mathcal{S}_t)$ is available for all the instances $t = 1, \ldots, T$ to train a machine learning model, $\text{ML}(\cdot)$, of our choice. In the next step, we take the parameters $\boldsymbol{\theta}_{\tilde{t}}$ of a new unseen problem instance $\tilde{t}$, and predict an invariant constraint set $\mathcal{S}_{\tilde{t}} = \text{ML}(\boldsymbol{\theta}_{\tilde{t}})$ with the already trained model. Finally, we run just a few iterations of the constraint generation method warm-started with the learned invariant constraint set $\mathcal{S}_{\tilde{t}}$. The prediction of the invariant constraint set and the constraint generation method are executed online. In contrast, the  strategy to build $\mathcal{S}_t$ for all training instances $t$ (Algorithm \ref{alg: identifying support constraints}), as well as the training of the machine learning algorithm $\text{ML}(\cdot)$ are performed offline. This way, the online computational burden is not affected. \textcolor{black}{Algorithm \ref{alg: pseudocode} shows a pseudocode of the main steps of our approach.}
\begin{algorithm}[!htb]
\captionsetup{labelfont={color=black, bf}}
\textcolor{black}{\\
\textbf{\underline{Offline phase:}}\\
\textbf{Input:} $\{(\bm{\theta}_t, \mathcal{B}_t)\},\, \forall t$.
\begin{enumerate}[label={\arabic*)}]
    \item For each train instance $t$:
    \begin{enumerate}
        \item Run Algorithm \ref{alg: identifying support constraints}.
        \item Obtain an invariant constraint set, $\mathcal{S}_t$.
    \end{enumerate}
    \item Train ML$(\cdot)$ using  $\{(\bm{\theta}_t, \mathcal{S}_t)\},\, \forall t$.
\end{enumerate}
\textbf{Output:} Trained ML methodology.\\
\textbf{\underline{Online phase:}}\\
\textbf{Input:} Trained ML strategy of previous step and $\bm{\theta}_{\tilde{t}}$ from a test instance $\tilde{t}$.
  \begin{enumerate}[start=3,label={\arabic*)}]
    \item Predict $\mathcal{S}_{\tilde{t}} = \text{ML}(\bm{\theta}_{\tilde{t}})$.
    \item Run CG initialized with the set $\mathcal{S}_{\tilde{t}}$.
     \end{enumerate}
\textbf{Output:} Optimal solution of the test problem instance $\tilde{t}$.
   \caption{\textcolor{black}{Pseudocode of the proposed methodology.}}
   \label{alg: pseudocode}
   }
\end{algorithm}

The main advantages of the proposed methodology are described below:

\begin{itemize}

\item[-] \textcolor{black}{Since we are running a constraint generation procedure that is warm-started with a carefully built constraint set, our approach retains the convergence optimality guarantees from the standard constraint generation method.}

%It retains the theoretical optimality guarantees from the standard constraint generation method since we are running a constraint generation procedure that is enhanced with an efficiently warm-started constraint set.

\item[-] The procedure to build the invariant constraint sets, $\mathcal{S}_t, \, \forall t$, of previously solved instances is also based on constraint generation. Hence, it is guaranteed to include all the non-binding constraints necessary to recover the optimal objective value. Therefore, all the past instances verify that the optimal solutions of the reduced problems are feasible and optimal for the original formulations.

\item[-] There is no condition about the machine learning algorithm that we apply in our approach. In other words, the sets of constraints $\mathcal{S}_t$ can be used for training \emph{any} machine learning method.

\item[-] The invariant constraint sets, $\mathcal{S}_t,\, \forall t$ and the machine learning algorithm, $\text{ML}(\cdot)$ are run offline. Therefore, the computational cost executed online to determine the optimal solution of a new unseen MILP instance is not affected.
\end{itemize}

\section{Computational Experiments} \label{sec: Computational Experience}
This section is devoted to the numerical experiments carried out in this paper. Section \ref{subs: Experimental setup} details the experimental setup, whereas Section \ref{subs: Case studies} explains the results derived from testing our proposal on two case studies.
  
\subsection{Experimental setup} \label{subs: Experimental setup}
  
To show the efficiency of our approach, we compare it with two algorithms. The first one is the standard constraint generation algorithm, denoted as CG,  which is based on pure optimization grounds and completely ignores the information provided by the data. In particular, for each instance, we sequentially add the violated constraints at each iteration. The second comparative approach is based on reference \cite{pineda2020data}. The authors propose a data-driven method where the constraint generation method is warm-started using only the information given by the binding constraint set, $\mathcal{B}$. We denote this method by $\mathcal{B}$-learner + CG. Finally, since we warm-start the constraint generation using an invariant constraint set $\mathcal{S}$, our methodology is denoted as $\mathcal{S}$-learner + CG. 
  
It is well-known that machine learning performance highly depends on the data division into training and test samples. Thus, to get stable out-of-sample results, leave-one-out is executed in this paper. More details about this technique can be found in \cite{hastie2009elements}. In particular, we assume given a database of $T$ previously solved MILP instances. As mentioned in Section \ref{sec: Invariant Constraint Set in MILPs}, we assume that the optimal solution of such instances is unique. Otherwise, retaining one of them is enough to our purposes. Note that, in the case of multiple solutions, solvers usually retain just one of them instead of the complete set of solutions. Hence, when multiple solutions appear, we only collect the one provided by the solver. 
  
The leave-one-out strategy consists in running $T$ iterations of our approach. At each iteration, we select one MILP instance and consider it as the test set, $\left\{\tilde{t}\right\}$ to be run in the online phase. The remaining $T-1$ instances constitute the training set, $\left\{1, \ldots, T\right\} \setminus\left\{\tilde{t}\right\}$. Then, we run the $\mathcal{S}$-learner + CG method. We first identify offline $T-1$ invariant constraint sets $\mathcal{S}_t$, for all training instance $t$ by running Algorithm \ref{alg: identifying support constraints}, and train (also offline) a machine learning model, $\text{ML}(\cdot)$, that learns an invariant constraint set in terms of the MILP parameters. The next step (to be also performed online) is to use the already trained model to predict an invariant constraint set for the test instance. Finally, such a predicted set is utilized to warm-start the constraint generation procedure that results in the optimal solution of the test instance. For the sake of comparison, both data-based strategies are applied in an equivalent way. Naturally, the $\mathcal{B}$-learner + CG procedure is trained with the series of $T-1$ binding constraints sets, $\mathcal{B}_t$, instead of the invariant constraint sets, $\mathcal{S}_t$. On the other hand, the pure optimization-based strategy CG does not need to divide the whole database into train and test sets. CG is run $T$ times, one per MILP instance. Each time, an iterative algorithm that sequentially adds the violated constraints is run. Figure \ref{fig: flowchart comparative algorithms} shows a scheme of our proposal \textcolor{black}{(given in Algorithm \ref{alg: pseudocode})} together with the two comparative algorithms, emphasizing the offline and online steps of each method.
  \begin{landscape}
 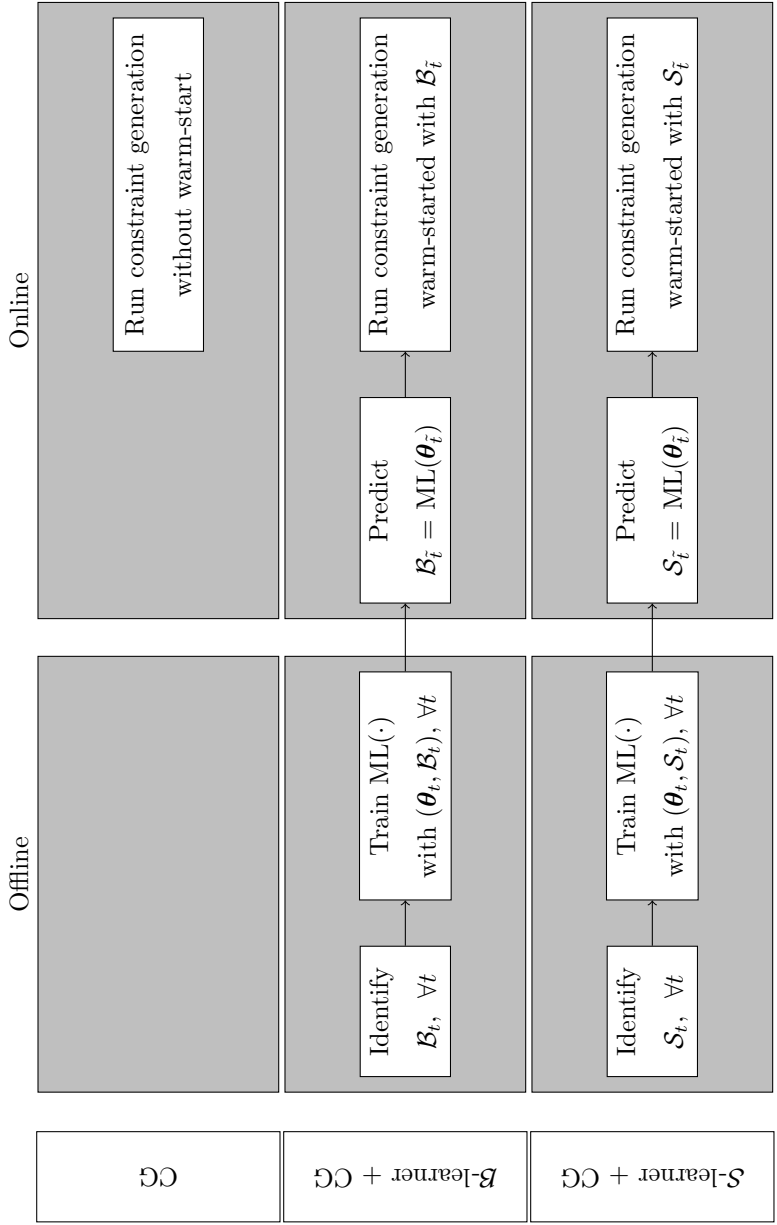
\begin{figure}[htb!]
 \hspace*{-1.5cm}
  %\centering
  %CG
  \begin{tikzpicture}
    \node [block, name=step1, text width=1.5cm, draw = none, fill=gray!50] {};
    \node [block, right of=step1, draw = none, fill=gray!50, text width=2.8cm ,node distance=3cm] (step2) {};
    \node [block, right of=step2,node distance=3.8cm, draw = none, fill=gray!50, text width=2.5cm] (step3) {};
    \node [block, right of=step3 , text width=4.2cm, node distance = 4.2cm] (step4) {Run constraint generation \\ without warm-start};
    \node[block, left of = step1, rotate = 90, node distance = 2.2cm] (CG){CG};
    \begin{scope}[on background layer]
    \node [container,fit=(step1) (step2), label = Offline] (offline) {};
    \node [container,fit=(step3) (step4), label = Online] (online) {};
\end{scope}
\end{tikzpicture}
%B-learner
\hspace*{-1.5cm}
\begin{tikzpicture}
    \node [block, name=step1, text width=1.5cm] {Identify \\ $\mathcal{B}_{t},\,\,\,
    \forall t$};
    \node [block, right of=step1, text width=2.8cm,node distance=3cm] (step2) {Train $\text{ML}(\cdot)$\\ with $(\bm{\theta}_t, \mathcal{B}_t), \, \forall t$};
    \node [block, right of=step2,node distance=3.8cm, text width=2.5cm] (step3) { Predict\\ $\mathcal{B}_{\tilde{t}}=\text{ML}(\bm{\theta}_{\tilde{t}})$};
    \node [block, right of=step3 , text width=4.2cm, node distance = 4.2cm] (step4) {Run constraint generation \\warm-started with $\mathcal{B}_{\tilde{t}}$};
    \node[block, left of = step1, rotate = 90, node distance = 2.2cm] (B-learner + CG){$\mathcal{B}$-learner + CG};
    \begin{scope}[on background layer]
    \node [container,fit=(step1) (step2)] (offline) {};
    \node [container,fit=(step3) (step4)] (online) {};
\end{scope}
    \draw [->] (step1) -- (step2);
    \draw [->] (step2) --  (step3);
    \draw [->] (step3) -- (step4);
\end{tikzpicture}
%S-learner
\hspace*{3cm}
  \begin{tikzpicture}
    \node [block, name=step1, text width=1.5cm] {Identify \\ $\mathcal{S}_{t},\,\,\,
    \forall t$};
    \node [block, right of=step1, text width=2.8cm, node distance=3cm] (step2) {Train $\text{ML}(\cdot)$\\ with $(\bm{\theta}_t, \mathcal{S}_t), \, \forall t$};
    \node [block, right of=step2, node distance=3.8cm, text width=2.5cm] (step3) { Predict\\ $\mathcal{S}_{\tilde{t}}=\text{ML}(\bm{\theta}_{\tilde{t}})$};
    \node [block, right of=step3 , text width=4.2cm, node distance = 4.2cm] (step4) {Run  constraint generation \\ warm-started with $\mathcal{S}_{\tilde{t}}$};
    \node[block, left of = step1, rotate = 90, node distance = 2.2cm] (S-learner + CG){$\mathcal{S}$-learner + CG};
    \begin{scope}[on background layer]
    \node [container,fit=(step1) (step2)] (offline) {};
    \node [container,fit=(step3) (step4)] (online) {};
\end{scope}
    \draw [->] (step1) -- (step2);
    \draw [->] (step2) --  (step3);
    \draw [->] (step3) -- (step4);
\end{tikzpicture}
  \caption{Flowchart of the three algorithms}
  \label{fig: flowchart comparative algorithms}
\end{figure}
\end{landscape}
The training of the machine learning model $\text{ML}(\cdot)$ and the subsequent construction of the set $\mathcal{S}_{\tilde{t}}$ (resp. $\mathcal{B}_{\tilde{t}}$)  is addressed as a binary classification problem. For this purpose, we assign a label $s_t^j = 1$ or $s_t^j = -1$ to each constraint $j\in\mathcal{J}$ of each problem instance $t$ with optimal solution $\boldsymbol{z}_{\bm{\theta}_t}^*[\mathcal{J}]$, depending on whether that constraint is in $\mathcal{S}_t$ (resp. $\mathcal{B}_{\tilde{t}}$)  or not, respectively. Likewise, each constraint $j\in\mathcal{J}$ in the new problem instance $\tilde{t}$ will be assigned the label $s_{\tilde{t}}^j = 1$ if the machine learning model $\text{ML}(\cdot)$ predicts that the constraint $j$ is in $\mathcal{S}_{\tilde t}$ (resp. $\mathcal{B}_{\tilde{t}}$), and label $s_{\tilde{t}}^j=-1$ otherwise.  Accordingly, any classification algorithm can, in principle, serve this purpose. In this paper, we select the well-known $k$ nearest neighbors, knn, due to its simplicity. See reference \cite{taunk2019brief} for more details in this regard. Nevertheless, alternative learning approaches such as Support Vector Machines, Neural Networks or Decision Trees can be applied as well. The work \cite{hastie2009elements} explains the main properties of these methodologies. 

Note that misclassifying has different consequences depending on the type of constraint we wrongly label. Indeed, adding a superfluous constraint into the set $\mathcal{S}_{\tilde{t}}$ is far much less damaging than failing to include a constraint in that set. The former case only leads to a slight increase in the size of the reduced problem, while the latter increases the risk that the solution to the reduced problem is infeasible in the original formulation, thus potentially increasing the number of iterations to be executed in the constraint generation method. For this reason, in this paper, we want to be on the conservative side when choosing the knn voting strategy. For a fixed $k$, we consider that a constraint $j$ of an unseen test instance, $\tilde{t}$, belongs to an invariant constraint set $\mathcal{S}_{\tilde{t}}$, i.e., $s^j_{\tilde{t}}=1$, if and only if at least one of the $k$ closest training instances includes such a constraint $j$ in the set $\mathcal{S}_t$. To ensure a fair comparison, we apply the same voting strategy when constructing $\mathcal{B}_{\tilde{t}}$ from $\mathcal{B}_{t},\,\forall t$ using the knn classification. In any case, a feature of the knn method is that the larger the $k$, the larger the size of $\mathcal{S}_{\tilde{t}}$  and $\mathcal{B}_{\tilde{t}}$. Thus, high values of $k$ are expected to result in predicted sets $\mathcal{S}_{\tilde{t}}$  and $\mathcal{B}_{\tilde{t}}$ with a larger number of constraints. Therefore, the number of iterations in the constraint generation method may be reduced at the expense of a potential increase in the running time of each iteration. \textcolor{black}{Then, there exists a trade-off between the value of $k$ and the possible increment of the computational burden. The user should decide the value of $k$ depending on their preferences.}
  
In general, any of the constraints of optimization problem \eqref{subeq: MILP} might be superfluous and thus unnecessary in a certain MILP instance. Nevertheless, it is often the case that, because of the nature and structure of the MILP under consideration, there is some group of constraints, say $\bar{\mathcal{J}} \subset \mathcal{J}$, which are more prone to be redundant and/or whose elimination from the original MILP  brings a substantial reduction in computational time. It may be, therefore, very useful to focus on that group of constraints only, and adding the set of constraints $\mathcal{J}\setminus\bar{\mathcal{J}}$ into $\mathcal{S}_{t}$ and $\mathcal{B}_{t}$, $\forall t$, by default. We notice that this is an issue analogous to that of deciding which group of binary variables is better to be learned in those strategies that help solve MILPs very fast by predicting the optimal value of some of these variables, see, for instance, \cite{lodi2020learning}. In the numerical experiments we present in Section \ref{subs: Case studies}, we will specify which constraints are considered in $\bar{\mathcal{J}}$.

\textcolor{black}{The efficiency of our proposal is measured using different performance metrics depending on the size of the datasets. For the toy example of Section \ref{subsubs: Toy Example}, we show the behaviour of our proposal $\mathcal{S}$-learner + CG compared with the alternative machine-learning-aided approach $\mathcal{B}$-learner + CG. We compare both strategies in terms of: i) the constraints used to warm-start the CG method, ii) the number of runs of CG, and iii) the final reduced set of constraints needed to get the optimal solution of the MILP. In constrast, for the large datasets of Section \ref{subsubs: Synthetic Example} and \ref{subsubs: Real-world Application: Unit Commitment Problem}, we measure the benefits of our approach} over the $T$ runs in terms of: i) the minimum and the maximum number of constraints considered in the reduced MILPs. Such a number of constraints is denoted as $C_{min}$, and $C_{max}$, respectively; ii) the minimum and maximum number of iterations executed in the constraint generation procedure of the three methodologies averaged over all test instances, denoted by $I_{min}$ and $I_{max}$, respectively; iii) the percentage, $P_1$,  of instances that require only one iteration of the constraint generation method; and iv) the percentage of online computational burden in comparison with the original MILP formulation, defined as \textcolor{black}{$\Delta = \frac{1}{T}\sum_{t=1}^T\delta_t$, where $\delta_t = 100\,\frac{\tau^{pred}_t+\tau^{CG}_t}{\tau^{MILP}_t}$}. For each instance $t$, $\tau^{MILP}_t$ denotes the time needed to solve the original MILP with the whole set of constraints $\mathcal{J}$. In addition,  $\tau^{pred}_t$ is the time  employed in predicting sets $\mathcal{B}_{\tilde{t}}$ or $\mathcal{S}_{\tilde{t}}$ by the knn strategy. Finally, the notation $\tau^{CG}_t$ indicates the computational time employed in all the iterations of the constraint generation method. 

\textcolor{black}{Apart from providing values that summarize the good performance of our approach, we have included three figures to illustrate the distribution of the performance measures. In particular, we provide the boxplots of the performance values obtained after running $T$ iterations for the different methodologies in terms of: i) the number of constraints considered in the reduced MILPs, ii) the number of iterations executed by the constraint generation strategy, and iii) the percentage of the computational time spent online compared to the original MILP formulation, $\delta_t$.}

All the experiments have been carried out on a cluster with $21$ Tb RAM, running Suse Leap 42 Linux distribution. MILP problems are coded in \texttt{Python 3.8} and \texttt{Pyomo 6.1.2} and solved using \texttt{Cplex 20.1.0}. \textcolor{black}{All the solver tuning parameters have been set to their default value except the \emph{mixed integer optimality gap tolerance} (\texttt{mipgap}) that has been fixed to $1e^{-10}$. Finally, in order to make our approach transparent, we have saved the data and the code of our proposal in \cite{OASYS2022warm}.}

%   \begin{figure}[htb!]
%  \hspace*{-0.7cm}
%   \centering
%   \begin{tikzpicture}
%     \node [block, name=input, text width=2.3cm] {\textbf{Input}: MILP\\ parameters, $\boldsymbol{\theta}$};
%     \node[block, name=B-learner, right of=input]{$\mathcal{B}$-learner};
%     \node[block, name=S-learner, below right=0.5cm and 0.6cm of input]{$\mathcal{S}$-learner};
%     \node[block, name=CG, right of=B-learner]{CG};
%     \node [block, right of=CG , text width=2cm] (output) {\textbf{Output}: optimal solution, $\boldsymbol{z^*_{\boldsymbol{\theta}}}$};

%     \draw [->] (input) |- (B-learner);
%     \draw [->] (B-learner) --  (CG);
%     \draw [->] (CG) -- (output);
%     \draw[->] ($(input)!.43!(B-learner)$) -- +(0, 1.5) -- + (5.5, 1.5) -| (CG.north);
%     %\draw[->] ($(input)!.43!(B-learner)$) .. controls +(up:2cm) and +(left:7mm).. (CG.north);
%     \draw[->] (input) -- ++(1.5cm,0) |- (S-learner);
%     \draw[->] (S-learner) -- ++(1.8cm,0) -| (CG);
% \end{tikzpicture}
%   \caption{Flowchart of the five strategies executed in the computational experiments.}
%   \label{fig: five algorithms numerical experience}
% \end{figure}

\subsection{Case Studies} \label{subs: Case studies}
The proposed methodology has been tested on \textcolor{black}{three} case studies: \textcolor{black}{A toy example (Section \ref{subsubs: Toy Example}), a} synthetic MILP (Section \ref{subsubs: Synthetic Example}) and a real-world application, the Unit Commitment Problem (Section \ref{subsubs: Real-world Application: Unit Commitment Problem}).

\subsubsection{\textcolor{black}{Toy Example}}\label{subsubs: Toy Example}
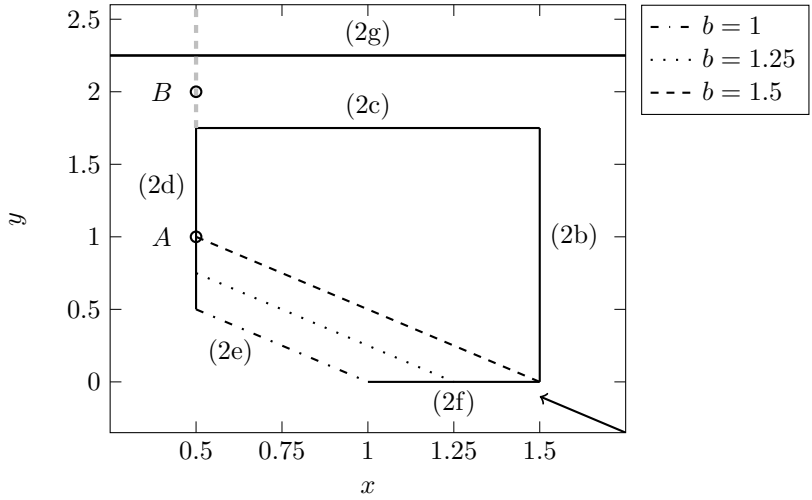
\begin{figure}
\captionsetup{labelfont={color=black}}
    \centering
    \begin{tikzpicture}[scale= 1,every node/.style={scale=1}]
	\begin{axis}[black, xmin=0.25,
	xmax=1.75,
	ymin=-0.35,
	ymax=2.6,
	xlabel=$x$,
	ylabel=$y$, xtick={0.5, 0.75, 1, 1.25, 1.5},
      xticklabels={0.5, 0.75, 1, 1.25, 1.5},
      ytick={0,  0.5,  1,  1.5,  2, 2.5},
      yticklabels={0,  0.5,  1,  1.5,  2, 2.5},
      legend pos=outer north east, 
      legend entries={$b = 1\phantom{.23}$, $b = 1.25$, $b = 1.5\phantom{3}$},
      legend style={draw=black}
              ]
  \pgfplotsset{tick style={black}}
  \addlegendimage{loosely dashdotted, thick, black},
  \addlegendimage{loosely dotted, thick, black},
  \addlegendimage{dashed, thick, black}]
    
	\addplot[thick, black] coordinates {(0.5,1.75)(1.5,1.75)};
	\addplot[thick, black] coordinates {(0.5,0.5)(0.5,1.75)};
	\addplot[thick, black] coordinates {(1.5,0)(1.5,1.75)};
	\addplot[thick, black] coordinates {(0.25,2.25)(1.75,2.25)};
	\addplot[thick, black] coordinates {(0.25,2.25)(1.75,2.25)};
	\addplot[thick, black] coordinates {(1,0)(1.5,0)};
	
	%\addplot[thick, lightgray, ultra thick] coordinates {(1,0)(1.5,0)};
	%\addplot[thick, lightgray, ultra thick] coordinates {(0.5,1)(1.5,1)};

	\addplot[thick,dashed, black] coordinates {(0.5,1)(1.5,0)};
	\addplot[thick,loosely dotted, black] coordinates {(0.5,0.75)(1.25,0)};
	\addplot[thick,loosely dashdotted, black] coordinates {(0.5,0.5)(1,0)};
	
	\addplot[dashed, lightgray, ultra thick] coordinates {(0.5,1.75)(0.5, 2.6)};
	%\addplot[dashed, blue!30, ultra thick] coordinates {(0.5,1.75)(0.5, 2.6)};

	\addplot[mark=o,thick, black] coordinates {(0.5,1)};
	\addplot[mark=o,thick, black] coordinates {(0.5,2)};
	
	\draw[->, thick, black] (axis cs: 1.75,-0.35) to (axis cs: 1.5,-0.1);

	\node[ black] at (axis cs: 0.4,1) {$A$};
	\node[ black] at (axis cs: 0.4,2) {$B$};
	\node[ black] at (axis cs: 1.6, 1) {(2b)};
	\node[ black] at (axis cs: 1, 1.9) {(2c)};
	\node[ black] at (axis cs: 0.4, 1.35) {(2d)};
	\node[ black] at (axis cs: 0.6, 0.2) {(2e)};
	\node[ black] at (axis cs: 1.25, -0.15) {(2f)};
	\node[ black] at (axis cs: 1, 2.4) {(2g)};
	
    \end{axis}
\end{tikzpicture}
    \caption{\textcolor{black}{Feasible region, direction of improvement of the objective function and optimal solution (point $A$) of the three instances obtained with the values $b\in\{1, 1.25,  1.5\}$ in Problem \eqref{subeq: toy example MILP}.}}
    \label{fig: toy example}
\end{figure}

\textcolor{black}{This section presents a toy example to illustrate how our proposal works.
To this aim, we formulate the optimization problem in \eqref{subeq: toy example MILP} with two decision variables, namely, $x\in\mathbb{R}$ and $y\in\mathbb{Z}$, and the feasible region given by the six constraints defined in \eqref{subeq: first constraint toy example MILP} - \eqref{subeq: sixth constraint toy example MILP}.}
\textcolor{black}{
\begin{subequations}\label{subeq: toy example MILP}
  \begin{empheq}[left=\empheqlbrace]{align}
    \min\limits_{\substack{{x\in\mathbb{R}},\\ y\in \mathbb{Z}}} & \, x-y\\
    \text{s.t. } & x\leq 1.5\label{subeq: first constraint toy example MILP}\\
    &y\leq 1.75 \label{subeq: second constraint toy example MILP}\\ 
    & x\geq 0.5\label{subeq: third constraint toy example MILP}\\
    & x+y\geq b \label{subeq: fourth constraint toy example MILP}\\
    & y\geq 0 \label{subeq: fifth constraint toy example MILP}\\ 
    & y\leq 2.25 \label{subeq: sixth constraint toy example MILP}
    \end{empheq}
  \end{subequations}
}
\textcolor{black}{
Note that problem \eqref{subeq: toy example MILP} depends on the parameter $b\in\mathbb{R}$ that appears in the right-hand side of constraint \eqref{subeq: fourth constraint toy example MILP}. Hence, in this toy example, we have $\bm{\theta} = b$. We assume given a database containing the results of three optimization problems solved for the values $b\in\{1, 1.25, 1.5\}$. Figure \ref{fig: toy example} shows the feasible region of the problem for the three different values of $b$. An arrow at the bottom right corner indicating the direction of improvement of the objective function is also depicted. It is easy to see that the optimal solution for these three problems is the point $A = (0.5, 1)$.}

\textcolor{black}{Now, let us assume that we have a new unobserved problem instance given by the value $b^{test} = 1.3$. The objective of this section is to learn a reduced set of constraints that will be used to initialize the constraint generation method and solve the problem for $b^{test}$. This way, the number of iterations needed to find the optimal solution of the MILP is reduced, and so does the associated computational burden. In this toy example, we illustrate that it is not enough to learn such a reduced set of constraints with the information given by the binding constraints set, $\mathcal{B}$.  We have to take into account also some crucial non-binding constraints included in the invariant constraint set $\mathcal{S}$ obtained after running Algorithm \ref{alg: identifying support constraints}. Table \ref{table: constraints included sets toy example} shows which are the constraints that belong to both sets $\mathcal{B}$ and $\mathcal{S}$ for the three training instances given by values $b\in\{1, 1.25, 1.5\}$. Note that, as stated in Section \ref{subs: Experimental setup}, the constraints belonging to the set $\mathcal{B}$ (resp. $\mathcal{S}$) are labeled with 1. On the other hand, the constraints that do not belong to $\mathcal{B}$ (resp. $\mathcal{S}$) are denoted with -1.}

\begin{table}[!htb]
\captionsetup{labelfont={color=black}}
        \centering
        \color{black}
        \begin{tabular}{cccccccc}
        \toprule
           &$b$&\eqref{subeq: first constraint toy example MILP}&\eqref{subeq: second constraint toy example MILP}&\eqref{subeq: third constraint toy example MILP} &\eqref{subeq: fourth constraint toy example MILP}&\eqref{subeq: fifth constraint toy example MILP}&\eqref{subeq: sixth constraint toy example MILP}\\
        \midrule
    \parbox[t]{2mm}{\multirow{3}{*}{$\mathcal{B}$}}&1&-1&-1&1&-1&-1&-1\\
    &1.25&-1&-1&1&-1&-1&-1\\
    &1.5&-1&-1&1&1&-1&-1\\
        \midrule
    \parbox[t]{2mm}{\multirow{3}{*}{$\mathcal{S}$}}&1&-1&1&1&-1&-1&-1\\
    &1.25&-1&1&1&-1&-1&-1\\
    &1.5&-1&1&1&1&-1&-1\\
             \bottomrule
        \end{tabular} \caption{\textcolor{black}{Constraints included in the sets $\mathcal{B}$ and $\mathcal{S}$ for the values $b\in\{1, 1.25, 1.5\}$ of the three training instances.}}\label{table: constraints included sets toy example}
    \end{table}
    
\textcolor{black}{The main difference between sets $\mathcal{B}$ and $\mathcal{S}$ is the constraint \eqref{subeq: second constraint toy example MILP},  which is included in the invariant constraint set $\mathcal{S}$ for the three values of $b$, but not in $\mathcal{B}$. Importantly, constraint \eqref{subeq: second constraint toy example MILP} is the unique non-binding constraint needed to recover the optimal solution given by the point $A$.  In effect, if this constraint is removed, then the optimal solution moves from point $A$ to point $B = (0.5, 2)$.}

\textcolor{black}{To find the initial set of constraints warm-starting the constraint generation method for the test problem instance given by $b^{test} = 1.3$, we train knn for $k\in\{1, 2, 3\}$ using the information provided by sets $\mathcal{B}$ and $\mathcal{S}$ in Table \ref{table: constraints included sets toy example}. Note that for $k = 1$, the closest problem instance to $b^{test}$ is the one associated to $b = 1.25$. In addition, for $k = 2$, the closest neighbors are the MILPs given by $b=1.25$ and $b = 1.5$. Finally, when $k = 3$, the three values of $b\in\{1, 1.25, 1.5\}$ are considered as the nearest neighbours of $b^{test}$. We collect in Table \ref{table: performance results toy example} the results obtained for both approaches in terms of: i) the constraints initially selected, ii) the number of iterations executed by the constraint generation strategy, and iii) the constraints employed to solve the reduced test MILP instance. Due to the small size of the problem, the computational burden of both approaches is negligible. Consequently, this information has been omitted.}

\begin{table}[!htb]
\captionsetup{labelfont={color=black}}
        \centering
        \color{black}
        \begin{tabular}{ccccc}
        \toprule
           &$k$&warm-start constraint set& iterations CG&final set of constraints\\
        \midrule
    \parbox[t]{2mm}{\multirow{3}{*}{\rotatebox[origin=c]{90}{\scriptsize$\mathcal{B}$-learner+CG}}}&1&\eqref{subeq: third constraint toy example MILP}&2&\eqref{subeq: second constraint toy example MILP}, \eqref{subeq: third constraint toy example MILP}\\
    &2&\eqref{subeq: third constraint toy example MILP}, \eqref{subeq: fourth constraint toy example MILP}&2&\eqref{subeq: second constraint toy example MILP}, \eqref{subeq: third constraint toy example MILP}, \eqref{subeq: fourth constraint toy example MILP}\\
    &3&\eqref{subeq: third constraint toy example MILP}, \eqref{subeq: fourth constraint toy example MILP}&2&\eqref{subeq: second constraint toy example MILP}, \eqref{subeq: third constraint toy example MILP}, \eqref{subeq: fourth constraint toy example MILP}\\
        \midrule
    \parbox[t]{2mm}{\multirow{3}{*}{\rotatebox[origin=c]{90}{\scriptsize$\mathcal{S}$-learner+CG}}}&1&\eqref{subeq: second constraint toy example MILP}, \eqref{subeq: third constraint toy example MILP}&1&\eqref{subeq: second constraint toy example MILP}, \eqref{subeq: third constraint toy example MILP}\\
    &2&\eqref{subeq: second constraint toy example MILP}, \eqref{subeq: third constraint toy example MILP}, \eqref{subeq: fourth constraint toy example MILP}&1&\eqref{subeq: second constraint toy example MILP}, \eqref{subeq: third constraint toy example MILP}, \eqref{subeq: fourth constraint toy example MILP}\\
    &3&\eqref{subeq: second constraint toy example MILP}, \eqref{subeq: third constraint toy example MILP}, \eqref{subeq: fourth constraint toy example MILP}&1&\eqref{subeq: second constraint toy example MILP}, \eqref{subeq: third constraint toy example MILP}, \eqref{subeq: fourth constraint toy example MILP}\\
             \bottomrule
        \end{tabular} \caption{\textcolor{black}{Performance results: Toy example}}\label{table: performance results toy example}
    \end{table}
    
\textcolor{black}{Several conclusions can be derived from the results shown in Table \ref{table: performance results toy example}. Since we are executing a modified version of the constraint generation method, it is guaranteed that the optimal solution of the problem instance given by $b^{test}$ is reached using both data-driven approaches. Indeed, such an optimal solution is also attained at point $A$ using a reduced set of two or three constraints (see the last column of Table \ref{table: performance results toy example}). This means that we have managed to decrease between 50\%-66\% the cardinality of the original set of constraints, depending on the value of $k$ we use. However, the number of iterations  that the CG method needs to perform to get this optimum varies depending on which of the two data-driven approaches is employed. Indeed, the number of CG iterations executed by the approach $\mathcal{S}$-learner + CG is smaller than those employed in $\mathcal{B}$-learner + CG for all values of $k$ (1 versus 2 iterations). This is due to the fact that the latter is trained just using the information taken from the binding constraints, which is not enough when there exist integer decision variables in the optimization problem. It is important to highlight that the approach $\mathcal{B}$-learner + CG is not able to find the optimal solution executing just one iteration of the CG method \emph{even if all the available training instances are used}, that is, even running the knn method with $k=3$. In contrast, our proposal $\mathcal{S}$-learner + CG always finds the optimal solution of the problem by only running the CG strategy once. This occurs because the non-binding constraint \eqref{subeq: second constraint toy example MILP} is included in the initial set of constraints used to warm-start the constraint generation method. Finally, as it was explained in Section \ref{subs: Experimental setup}, we want to be on the conservative side when choosing the knn voting choice. This is the reason why constraint \eqref{subeq: fourth constraint toy example MILP} is included in the initial set of constraints of both approaches for $k = 2$ and $k = 3$, even if such a constraint is not necessary to reach the optimal solution of the test problem instance. However, including such a small number of non-critical constraints have minor consequences in terms of computational times, as will be observed in the larger datasets of Sections \ref{subsubs: Synthetic Example} and \ref{subsubs: Real-world Application: Unit Commitment Problem}.}

% \textcolor{blue}{As an illustration of our proposal, this section shows .... \\
% - Toy example to illustrate our approach.\\
% - Optimization problem with two variables and six constraints.\\
% - Constraint 4 depends on one parameter.\\
% - Figure 2 shows the feasible region, direction of optimization, optimum value.\\
% - Note that constraints XXX are binding, and constraints XXX non-binding but crucial. Otherwise, optimum goes to B. Distinguish the different values of b.\\
% - If a new test instance come with b = 1.1 or b = 1.4.\\
% - Include table of results for the labels of the different constraints.\\
% - We run our approach for k = 3 (most conservative approach).\\
% - Note that when only one constraint is binding, then it is selected.\\
% - In the reduced MILP, then we predict the set S as XXX.\\
% - In this particular case, the number of iterations is 1.\\
% - Computational time makes no sense here.
%}

\subsubsection{Synthetic Setup} \label{subsubs: Synthetic Example}
    
In this section, we restrict ourselves to the MILP problem of the form \eqref{subeq: binary MILP}:
    
\begin{subequations}\label{subeq: binary MILP}
  \begin{empheq}[left=\empheqlbrace]{align}
    \min\limits_{{\boldsymbol{x}\in\mathbb{R}^n},\, \boldsymbol{y}\in\{0,\,1\}^{n}}&\sum\limits_{i = 1}^nc_ix_i\\
    \text{s.t. }  & \sum\limits_{i= 1}^n a_{ij}x_i\leq b_j, \quad j =1, \ldots, m \label{subeq: inequality constraint}\\ 
    &l_iy_i \leq x_i\leq u_iy_i,\quad i=1, \ldots, n \label{subeq: box constraints}
    \end{empheq}
  \end{subequations}
  
\noindent where $\boldsymbol{a}_i = (a_{i1}, \ldots, a_{im})^{\intercal}, \forall i \leq n$, and $\boldsymbol{b} = (b_1, \ldots, b_m)^{\intercal}$ are column vectors in $\mathbb{R}^{m}$, and $\boldsymbol{c} = (c_1, \ldots, c_n)^{\intercal},\, \boldsymbol{l}= (l_1, \ldots, l_n)^{\intercal}$ and $\boldsymbol{u} =(u_1, \ldots, u_n)^{\intercal}$ are column vectors in $\mathbb{R}^n$. 

MILPs like \eqref{subeq: binary MILP} can be interpreted as linear programs where some of the continuous variables $x_i$ have a forbidden zone within the range $(0, l_i)$. Consequently, problems like \eqref{subeq: binary MILP} contain the so-called logical constraints, where a continuous variable vanishes if the associated binary variable is zero. This type of problems, with a logical relationship between continuous and binary variables, has a wide variety of applications, as the authors of \cite{bertsimas2021unified} explain. For a real-life example, one can think in a nuclear energy context. For instance, a nuclear unit whose maximum power is $1000$ MW cannot generate energy within the range $[0, 500]$ due to the nuclear reactor stability.
 
As mentioned in Section \ref{subs: Experimental setup}, it may be computationally productive to screen out only a subset of constraints $\bar{\mathcal{J}} \subset \mathcal{J}$. In the case of Problem~\eqref{subeq: binary MILP}, for example, most of the constraints \eqref{subeq: box constraints} are expected to be binding at the optimal solution. Therefore, we consider in this example that $\bar{\mathcal{J}}$ is solely formed by the $m$ constraints in \eqref{subeq: inequality constraint}. Consequently, the $n$ constraints in \eqref{subeq: box constraints} are included by default into  $\mathcal{S}_{t}$ and $\mathcal{B}_{t}$, $\forall t$.

We assume given a database with $T = 1000$ optimization problems of the type of \eqref{subeq: binary MILP}. Each optimization problem $t$ comprises $m = 250$ constraints, $n=500$ continuous variables, and $n=500$ binaries. We assume that the problems just depend on the parameter $\bm{\theta} = \boldsymbol{b}$, i.e., $\boldsymbol{a}_i, \forall i$, $\boldsymbol{c}$, $\boldsymbol{l}$, and $\boldsymbol{u}$ remain fixed for the $1000$ optimization problems.
   
To synthetically generate the database, the values for the parameters $\boldsymbol{a}_i, \forall i$, $\boldsymbol{c}$, $\boldsymbol{l}$, and $\boldsymbol{u}$ (which the 1000 MILPs share) have been randomly selected according to a normal distribution with mean $0$ and standard deviation $10$, i.e., $\mathcal{N}(0,10)$. Note that we assure that lower bounds $l_i$ take on smaller values than the upper bounds $u_i$, i.e.,  $l_i< u_i,\, i = 1, \ldots, n$. Then, $1000$ parameter vectors $\boldsymbol{b}$ have been generated again according to the same distribution $\mathcal{N}(0,10)$. The entire database can be downloaded from \cite{OASYS2022warm}.

\begin{figure}
\captionsetup{labelfont={color=black}}
\thisfloatpagestyle{empty}
\vspace*{-3cm}
\begin{subfigure}{\textwidth}
\centering
    \includegraphics[scale=0.6]{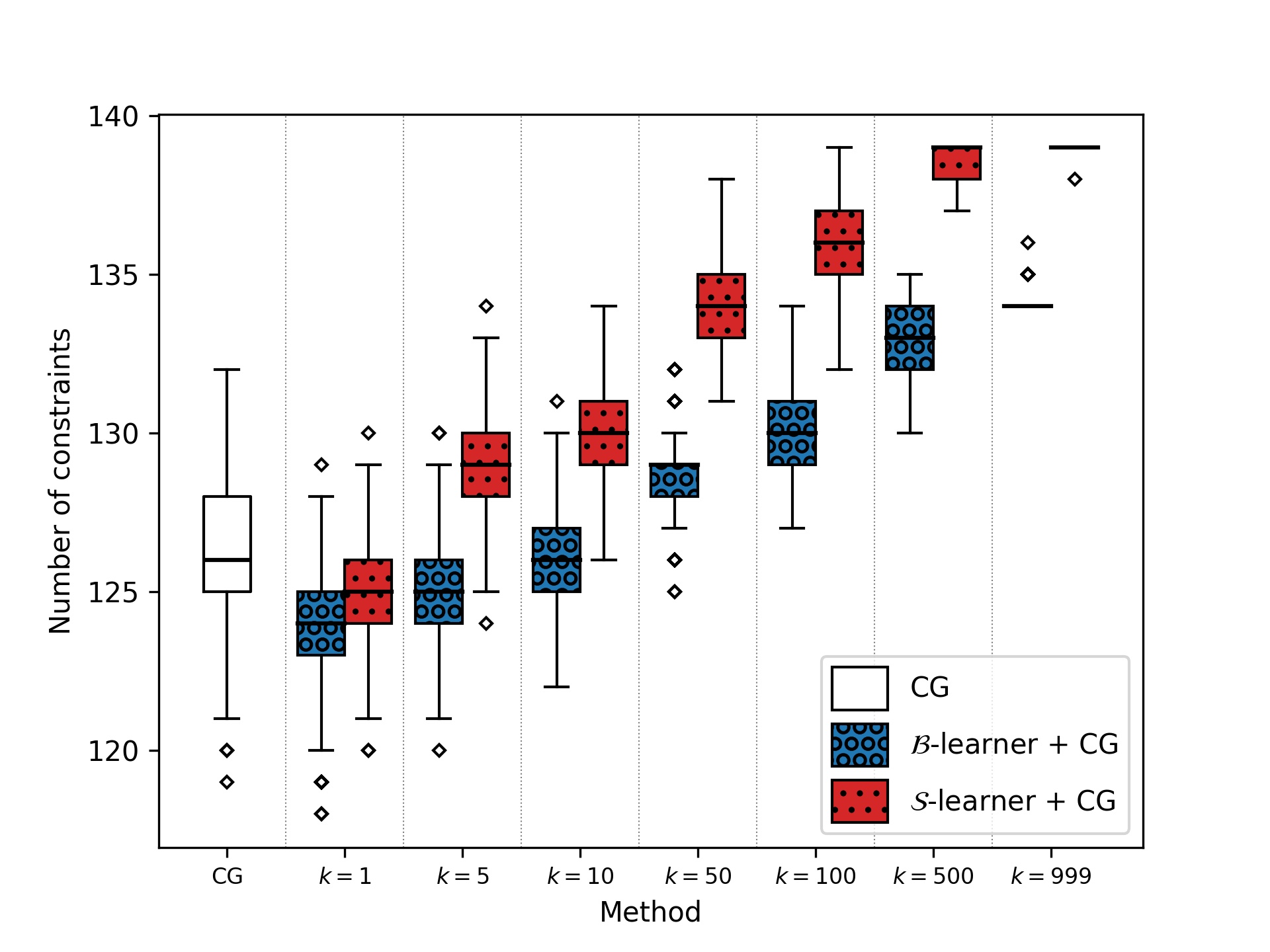}
    \caption{\textcolor{black}{Number of constraints considered in the reduced MILPs.}}
    \label{subfigure: number of constraints MILP}
\end{subfigure}
\\
\begin{subfigure}{\textwidth}
\centering
    \includegraphics[scale=0.6]{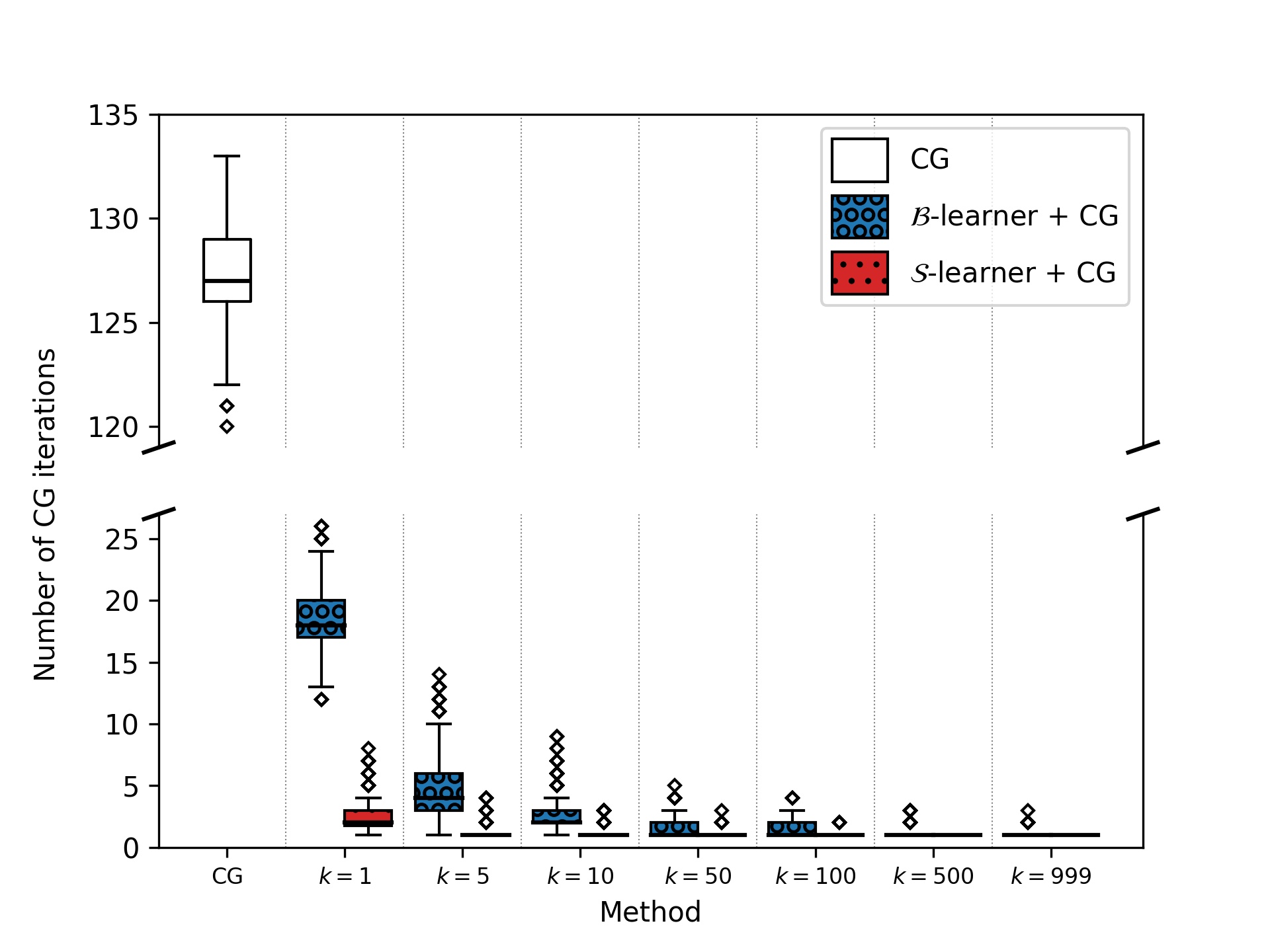}
    \caption{\textcolor{black}{Number of CG iterations.}}
    \label{subfigure: number of CG iterations MILP}
\end{subfigure}
\\
\begin{subfigure}{\textwidth}
\centering
    \includegraphics[scale=0.6]{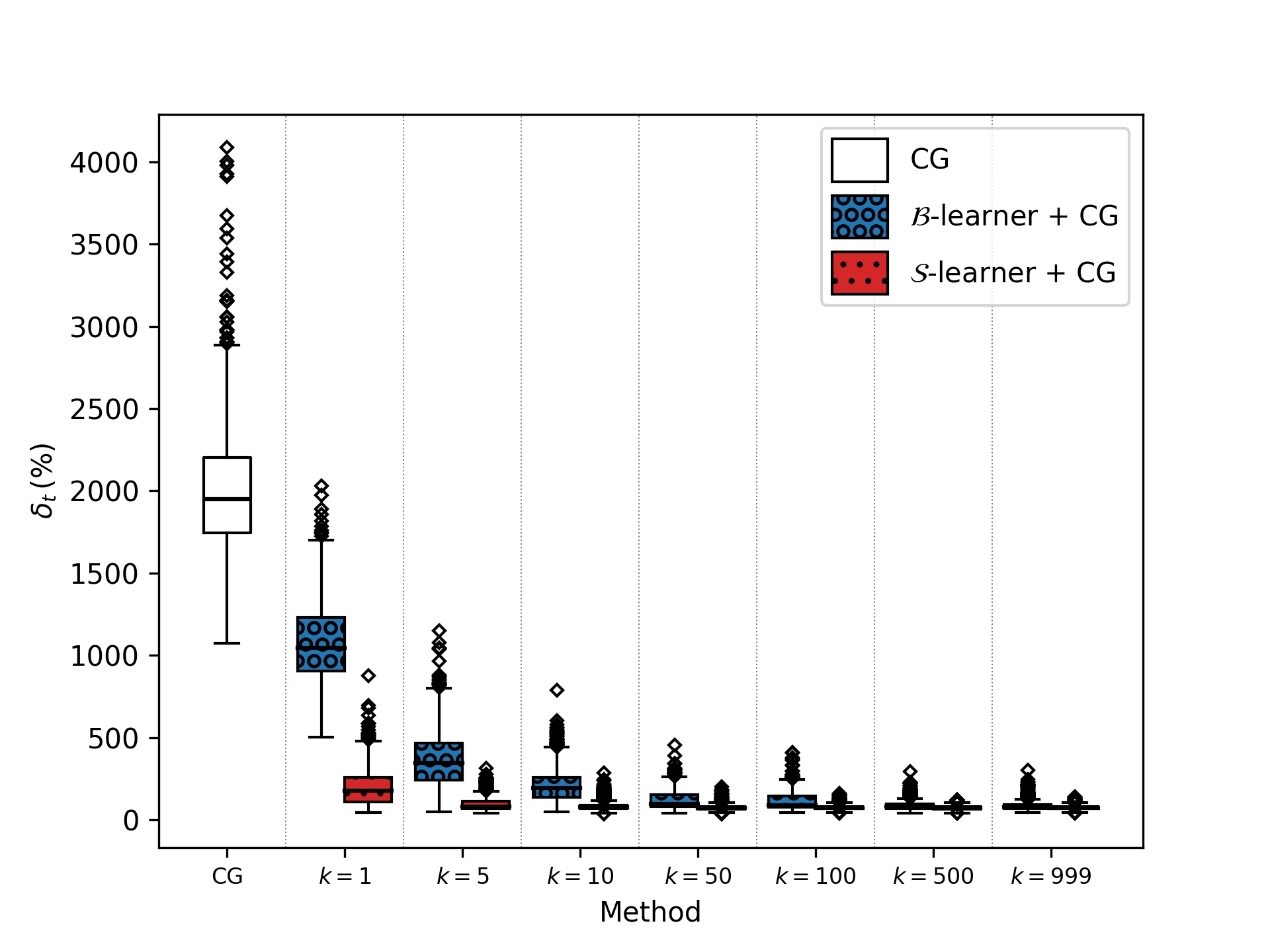}
    \caption{\textcolor{black}{Percentage of online computational burden in comparison with the original MILP formulation.}}
    \label{subfigure: online time MILP}
\end{subfigure}
\caption{\textcolor{black}{Performance results: Synthetic MILP.}}
\label{figure: performance results MILP}
\end{figure}

\textcolor{black}{Figure~\ref{figure: performance results MILP} and }Table~\ref{table: performance results MILP} shows the performance metrics for the three methodologies. The data-driven strategies $\mathcal{B}$-learner + CG and $\mathcal{S}$-learner + CG include output results after running the knn algorithm for different values of $k\in\{1, 5, 10, 50, 100,$ $500, 999\}$. Training the knn algorithm with $k = 999$ is equivalent to running a naive method that includes a constraint $j$ in the set $\mathcal{S}_{\tilde{t}}$  (resp. $\mathcal{B}_{\tilde{t}}$) if that constraint is contained in at least one of the training sets $\mathcal{S}_{t}, \, \forall t$  (resp. $\mathcal{B}_{t}, \, \forall t$). Note also that due to their nature, it is obvious that the three algorithms recover the optimal solution of the original MILP instances. %In particular, our proposal $\mathcal{S}$-learner + CG guarantees to attain the optimal solution of the problem.

\begin{table}[!htb]
        \centering
        \begin{tabular}{cccccc}
        \toprule
           &$k$&$[C_{min}, C_{max}]$&$[I_{min}, I_{max}]$ &$P_1 (\%)$&$\Delta (\%)$\\
          \midrule 
          \parbox[t]{2mm}{\multirow{1}{*}{\rotatebox[origin=c]{90}{CG}}} &-& [119, 132] &[120, 133]&0.0&1956.20\\
        \midrule
    \parbox[t]{2mm}{\multirow{7}{*}{\rotatebox[origin=c]{90}{$\mathcal{B}$-learner+CG}}}&1&[118, 129]&[12, 26]&0.0&1050.89\\
    &5&[120, 130]&[1, 14]&2.8&366.71\\
    &10&[122, 131]&[1, 9]&20.2&204.87\\
    &50&[125, 132]&[1, 5]&58.4&119.76\\
    &100&[127, 134]&[1, 4]&66.4&112.33\\
    &500&[130, 135]&[1, 3]&83.8&87.90\\
    &999&[134, 136]&[1, 3]&85.6&87.55\\
        \midrule
        \parbox[t]{2mm}{\multirow{7}{*}{\rotatebox[origin=c]{90}{$\mathcal{S}$-learner+CG}}}&1&[120, 130]&[1, 8]&25.0&195.67\\
    &5&[124, 134]&[1, 4]&77.7&95.25\\
    &10&[126, 134]&[1, 3]&90.6&81.06\\
    &50&[131, 138]&[1, 3]&98.3&73.25\\
    &100&[132, 139]&[1, 2]&99.3&74.47\\
    &500&[137, 139]&[1, 1]&100.0&72.22\\
    &999&[138, 139]&[1, 1]&100.0&74.38\\
             \bottomrule
        \end{tabular} \caption{Performance results: Synthetic MILP.}\label{table: performance results MILP}
    \end{table}

It can be observed that the number of iterations that need to be executed by the optimization-based method CG is significantly larger than those required by the data-driven methods $\mathcal{B}$-learner + CG and $\mathcal{S}$-learner + CG. \textcolor{black}{Indeed, the $y$-axis of Figure~\ref{subfigure: number of CG iterations MILP} has been divided into two different parts. This way, the large number of iterations needed in the pure optimization-based method CG does not affect the visualization of the number of CG iterations in the data-aided approaches. In addition, such a difference in the number of iterations} clearly shows the benefits of using machine learning tools to alleviate the online computational burden in contrast with the use of pure optimization-assisted strategies.

Moreover, if we compare the number of CG iterations in both data-driven methods, $\mathcal{B}$-learner + CG and $\mathcal{S}$-learner + CG, we can see that the number of iterations to be run in the former approach is larger than in the latter, for all the values of $k$. \textcolor{black}{Particularly, if we focus on the value $k = 1$ of Figure~\ref{subfigure: number of CG iterations MILP}, we observe that the minimum number of iterations that the method $\mathcal{B}$-learner + CG executed is larger than the maximum number of iterations of the strategy $\mathcal{S}$-learner + CG. Indeed, we can affirm that} in order to reduce the number of iterations of the constraint generation method, it is important to find a good initial set of constraints for the warm-start. Actually, if such an initial set is not good enough, solving the original MILP to optimality with the whole set of constraints, $\mathcal{J}$, could be better in terms of the computational burden than running a few iterations of the reduced MILP instances. For example, observe the output results \textcolor{black}{of Table~\ref{table: performance results MILP}} of $k = 100$ for both approaches $\mathcal{B}$-learner + CG and $\mathcal{S}$-learner + CG. It can be seen that, in the first case, the computational load is increased around $12\%$ \textcolor{black}{on average}, whereas in the second case, a reduction of approximately $25$ percentage points is attained. In this regard, we should also highlight \textcolor{black}{from Table~\ref{table: performance results MILP} } that running $\mathcal{B}$-learner + CG with values $k\in\{1, 5, 10, 50, 100\}$ provides no online running time benefits. In contrast, online computational savings can be observed from $k = 5$ if $\mathcal{S}$-learner + CG is performed. That is, the use of the sets $\mathcal{S}_t,\, \forall t$ which are built offline, is advantageous from a computational point of view.

In addition, we emphasize the results \textcolor{black}{from Table~\ref{table: performance results MILP}} of $k=500$ and $k=999$ for the $\mathcal{S}$-learner + CG. Note that the number of CG iterations in both cases is exactly one. That means that the prediction of the invariant constraints sets includes all the binding and non-binding constraints necessary to recover the optimal solution of the original MILP formulation. It is important to remark that the $\mathcal{B}$-learner + CG is unable to reproduce these results even for $k=999$, i.e., even using a naive method with all the available data. The reason for this issue is that the strategy $\mathcal{B}$-learner + CG is only based on the binding constraints sets, $\mathcal{B}_t, \forall t$, which are not sufficient to recover the optimal solution. Therefore, the risk of generating reduced MILPs, which are not equivalent to the original ones, increases.

Regarding the number of constraints of the reduced MILP formulations, it can be observed \textcolor{black}{in Table~\ref{table: performance results MILP} and Figure~\ref{subfigure: number of constraints MILP}} that the number of constraints necessary for solving the unseen instances includes around $50\%$-$55\%$ of the total number of constraints, independently of which of the three algorithms is run. This means that the size of the original optimization problem is considerably reduced with any of the methodologies. Regarding the data-driven methods, the number of retained constraints in the  $\mathcal{S}$-based methods is, as expected, slightly greater than in the $\mathcal{B}$-based ones. The difference boils down to a few extra constraints, which, thus, barely increases the size of the reduced MILPs. Very importantly, however, leaving these few constraints in the reduced MILPs has a major impact on the number of CG iterations, and thus in the online computational load, as Table \ref{table: performance results MILP} \textcolor{black}{and Figure~\ref{subfigure: online time MILP}} show.

In this vein, it is essential to remark that the simple fact of including any type of constraint in the set $\mathcal{S}_t$ is not enough to recover the optimal solution in a shorter computational time. For instance, notice \textcolor{black}{in Figure~\ref{subfigure: number of constraints MILP}} that the \textcolor{black}{distribution} of retained constraints after training the $\mathcal{B}$-learner + CG with $k = 100$, on the one hand, and the $\mathcal{S}$-learner + CG with $k = 10$, on the other, is very similar. However, \textcolor{black}{looking at Table~\ref{table: performance results MILP} we can see that} there is no online time reduction in the first case, whereas around $80\%$ of the computational load is employed in the second case. This is because more critical non-binding constraints are retained in the reduced MILPs in the latter case.

To sum up, our strategy is able to attain the optimal solution of an MILP, thanks to the efficient warm-start of the constraint generation method. Actually, we are able to correctly identify offline an invariant constraint set of an MILP, $\mathcal{S}_t,\, \forall t$. In doing so, the performance of the machine learning tool (knn in our case) that is used to initialize the constraint generation method is improved. This improvement substantially decreases the number of iterations performed at the expense of a slight increment in the cardinality of the set of retained constraints. This increment, however, does not involve an increase in the online solution time. In addition, the resulting reduced problems are easier to solve than the original MILPs.

\subsubsection{Real-world Application: Unit Commitment Problem}\label{subsubs: Real-world Application: Unit Commitment Problem}
    
The Unit Commitment problem (UC) is one of the most important problems in power systems, as the authors in  \cite{gomez2018electric} affirm. The goal of UC is to determine, at minimum cost, the on/off status and the power to be dispatched by each generation unit in order to satisfy the electric demand. Mathematically, the (DC version of the) UC problem can be formulated as the following MILP:
\begin{subequations}\label{subeq: UC}
  \begin{empheq}[left=\empheqlbrace]{align}
    \min\limits_{\boldsymbol{x}\in\mathbb{R}^{n},\, \boldsymbol{y}\in\{0,\, 1\}^{n}}&\sum\limits_{i=1}^n c_ix_i \label{subeq: cost function UC}\\
    \text{s.t. }  & \sum\limits_{i=1}^nx_i= \sum\limits_{i = 1}^n d_i, \label{subeq: balance constraint}\\
    &-f_j\leq \sum\limits_{i = 1}^n a_{ij}(x_i-d_i)\leq f_j, \quad  j = 1, \ldots, m\label{subeq: line constraints}\\
    &l_iy_i \leq x_i\leq u_iy_i,\quad i=1, \ldots, n \label{subeq: bound constraints power}
    \end{empheq}
  \end{subequations}

\noindent where $x_i$ is the power dispatched of generator $i$ and $y_i$ is a binary variable indicating whether the generator is turned on or turned off. In addition, $c_i$ is the marginal cost of generator $i$, $d_i$ is the electric demand at node $i$, $a_{ij}$ are the so-called Power Transfer Distribution Factors (PTDF) in \cite{gomez2018electric}, and $f_j$, $l_i$ and $u_i$ are, respectively, flow and power generation limits. The objective function \eqref{subeq: cost function UC} aims to minimize the total cost. Constraint \eqref{subeq: balance constraint} is the power balance equation, enforcing the supply of the total demand in the power network. Constraints \eqref{subeq: line constraints} limit the flow of line $j$, given by $\sum_{i = 1}^n a_{ij}(x_i-d_i)$, within the range $[-f_j,\,f_j]$. Finally, constraints \eqref{subeq: bound constraints power} ensures that the power dispatched $x_i$  be within $l_i$ and $u_i$ if and only if generator $i$ is turned on, i.e. if and only if $y_i=1$. We remark that, for simplicity, formulation~\eqref{subeq: UC} considers that there is at most one generator connected to each network node.
  
The Unit Commitment problem is a suitable application to test the performance of our method for three reasons. First, the increasing integration of renewable sources in current power systems requires that the unit commitment problem be solved multiple times within short-time windows so that commitment decisions can be adapted to rapid changes in operating conditions. In practice, this means that this problem must be solved as fast as possible. Second, the unit commitment problem is solved several times per day with only minor changes in the input data. Therefore, historical data of previous instances are usually available to be used in learning tasks. Third, implementing commitment decisions that violate some of the security constraints of the UC may lead to catastrophic events such as power blackouts. Therefore, attaining optimal (and feasible) solutions is also a requirement for this practical application.

While the marginal cost and production limits of power plants and network topology do not typically change over a year, the electric demand suffers from daily and weekly fluctuations. Hence, we decide to fix parameters $a_{ij}, l_i, u_i, f_j$ and $c_i$, $\forall i, j$, and just vary the input parameter  $\bm{\theta} = \boldsymbol{d}$. Moreover, it is known that, in practice, just a small percentage of the power flow constraints  \eqref{subeq: line constraints} are binding at the optimum for typical power systems. See reference \cite{bouffard2005umbrella} for further details. Hence, we consider in this example that $\bar{\mathcal{J}}$ is made up of the $2m$ constraints in \eqref{subeq: line constraints}. The entire dataset consists of $n=96$ continuous and $n = 96$ binary variables and $m =120$, leading to a total of $240$ constraints \eqref{subeq: line constraints}. The number of problem instances is $T = 8640$, corresponding to $360$ days of data measured every hour. More details about this data can be found in \cite{OASYS2022warm}.
 
\begin{figure}
\captionsetup{labelfont={color=black}}
\thisfloatpagestyle{empty}
\vspace*{-3cm}
\begin{subfigure}{\textwidth}
\centering
    \includegraphics[scale=0.6]{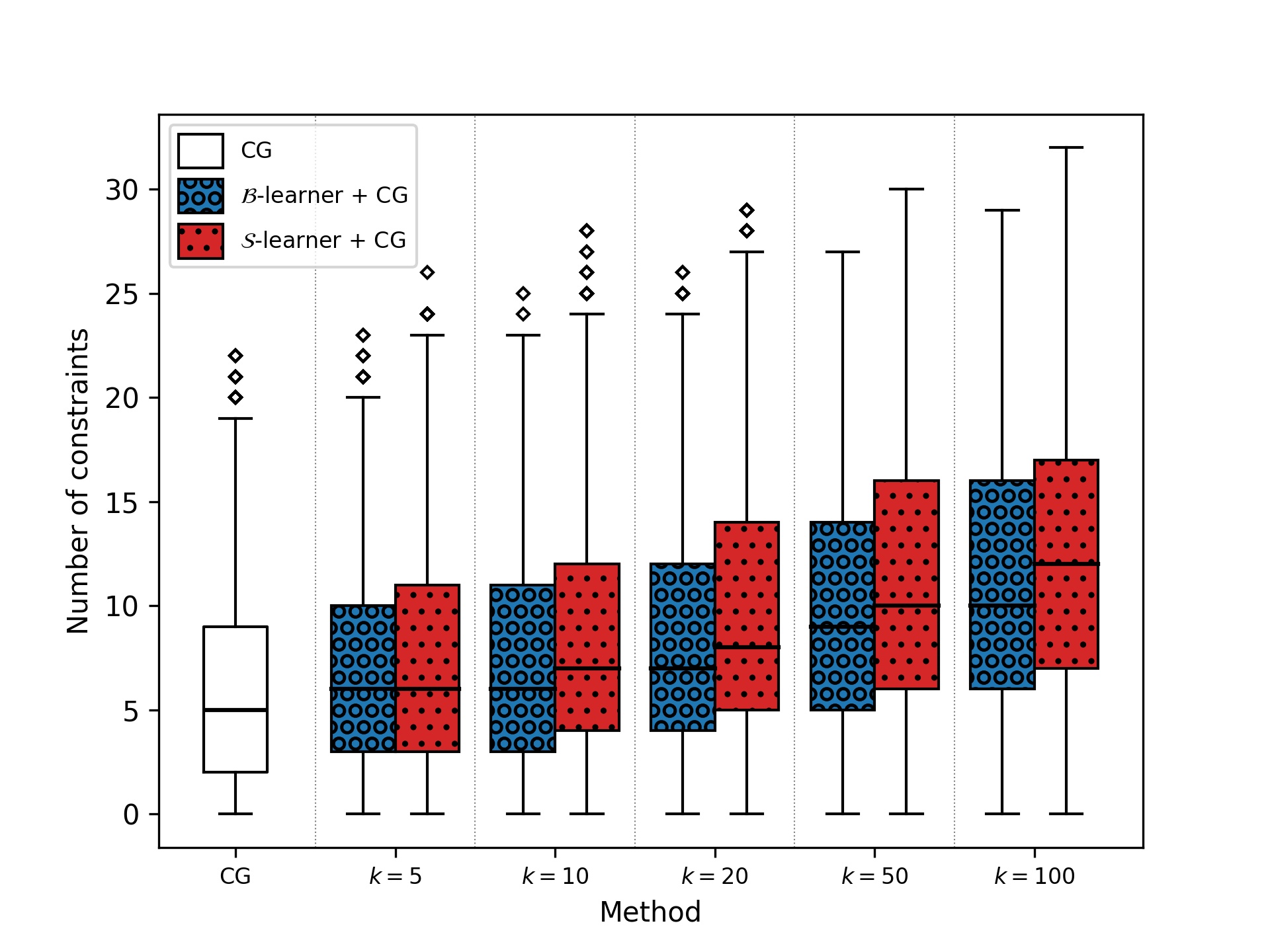}
    \caption{\textcolor{black}{Number of constraints considered in the reduced MILPs.}}
    \label{subfig: number of constraints UC}
\end{subfigure}
\\
\begin{subfigure}{\textwidth}
\centering
    \includegraphics[scale=0.6]{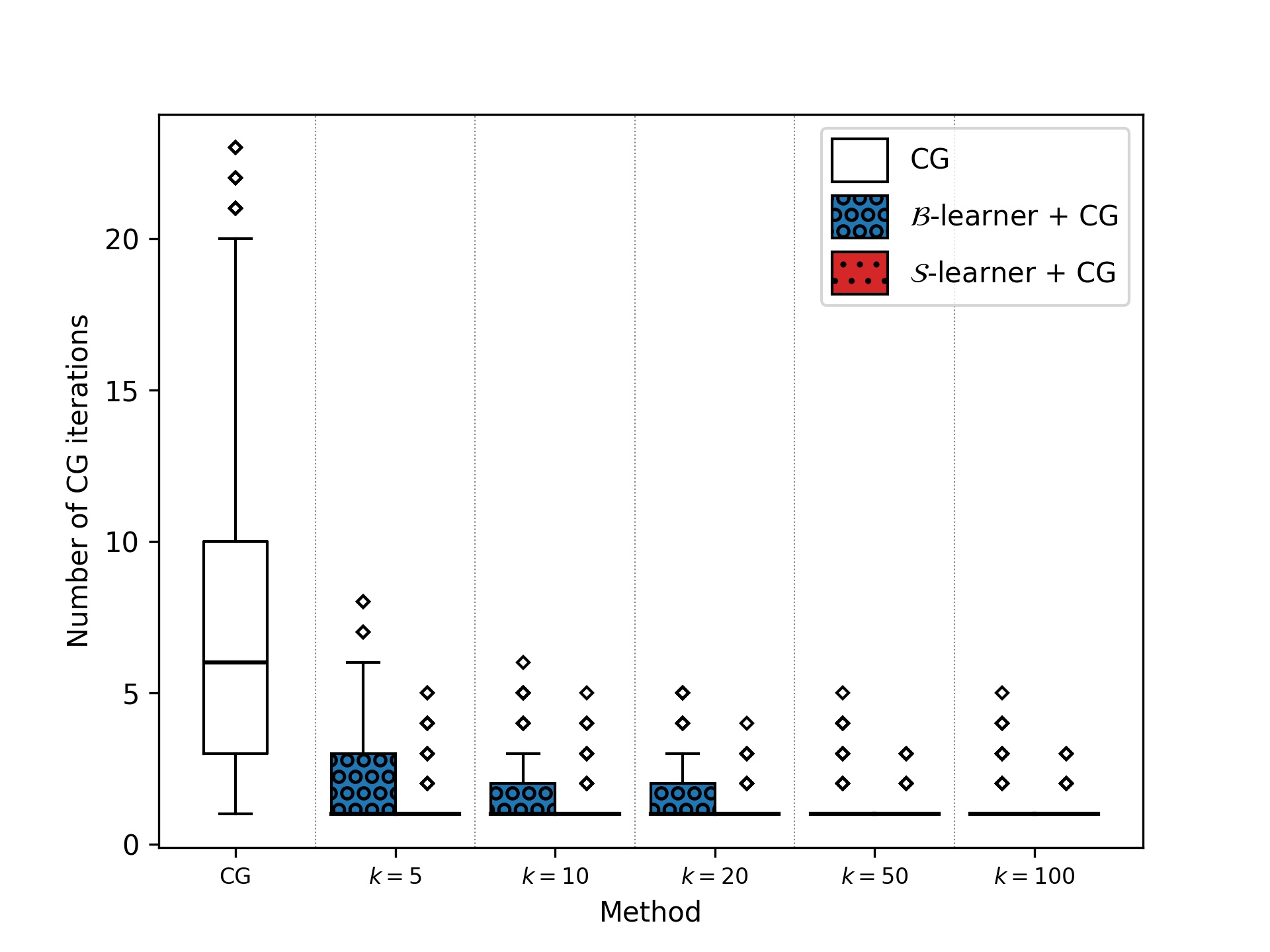}
    \caption{\textcolor{black}{Number of CG iterations.}}
    \label{subfig: number of CG iterations UC}
\end{subfigure}
\\
\begin{subfigure}{\textwidth}
\centering
    \includegraphics[scale=0.6]{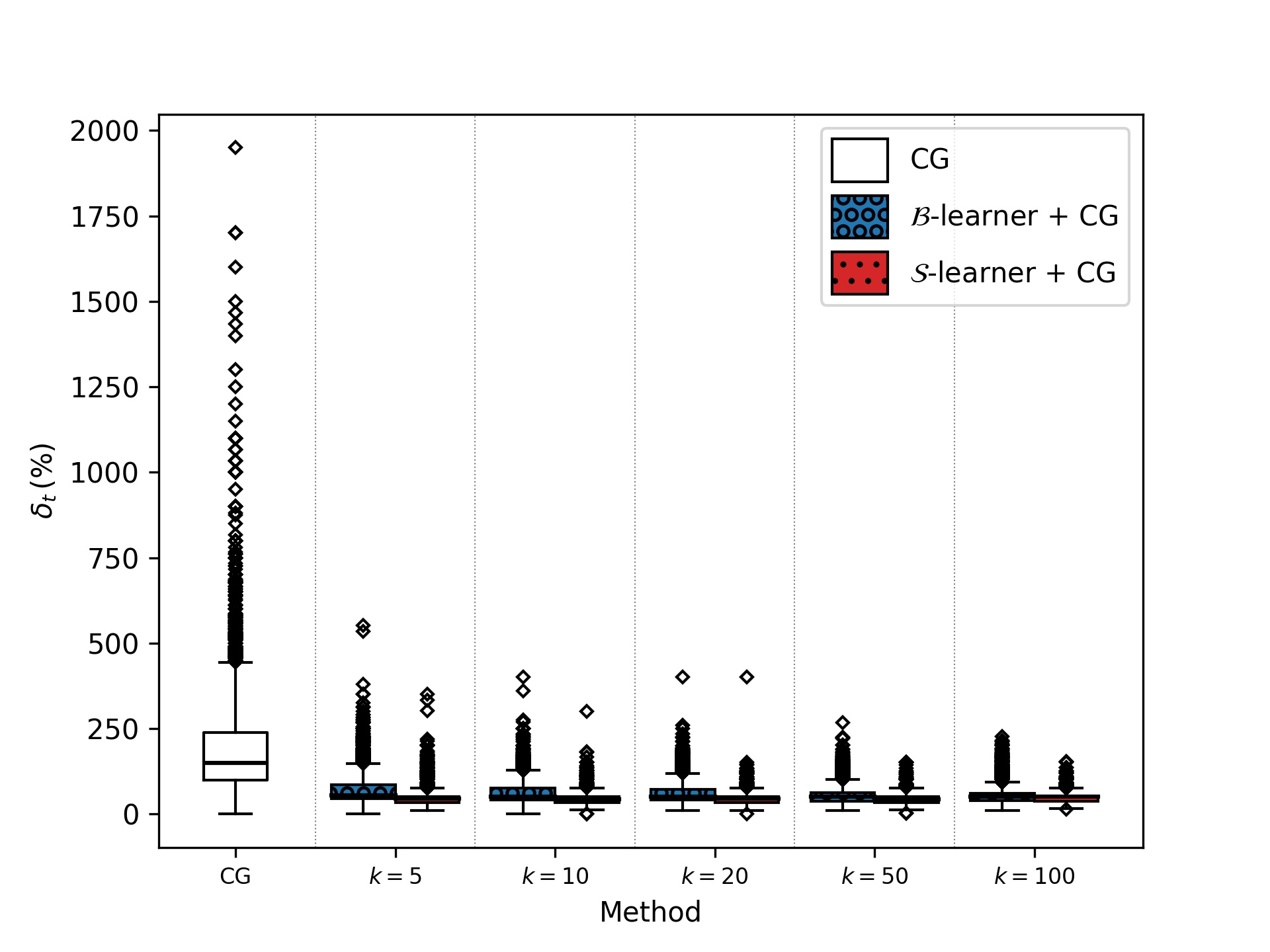}
    \caption{\textcolor{black}{Percentage of online computational burden in comparison with the original MILP formulation.}}
    \label{subfig: online time UC}
\end{subfigure}
\caption{\textcolor{black}{Performance results: Unit Commitment.}}
\label{fig: performance results UC}
\end{figure}

The performance metrics of our methodology are collated in Table \ref{table: performance results UC} \textcolor{black}{and Figure~\ref{fig: performance results UC}} for $k\in\{5, 10, 20, 50, 100\}$.

\begin{table}[!htb]
        \centering
        \begin{tabular}{cccccc}
        \toprule
           &$k$&$[C_{min}, C_{max}]$&$[I_{min}, I_{max}]$ &$P_1 (\%)$&$\Delta (\%)$\\
          \midrule 
          \parbox[t]{2mm}{\multirow{1}{*}{\rotatebox[origin=c]{90}{CG}}}&-& [0, 22] &[1, 23]&9.16&188.38\\
        \midrule
    \parbox[t]{2mm}{\multirow{5}{*}{\rotatebox[origin=c]{90}{$\mathcal{B}$-learner+CG}}}&5&[0, 23]&[1, 8]&54.40&74.09\\
    &10&[0, 25]&[1, 6]&62.01&63.19\\
    &20&[0, 26]&[1, 5]&68.28&62.35\\
    &50&[0, 27]&[1, 5]&76.90&54.56\\
    &100&[0, 29]&[1, 5]&83.70&54.62\\
        \midrule
    \parbox[t]{2mm}{\multirow{5}{*}{\rotatebox[origin=c]{90}{$\mathcal{S}$-learner+CG}}}&5&[0, 26]&[1, 5]&92.66&44.84\\
    &10&[0, 28]&[1, 5]&97.21&40.57\\
    &20&[0, 29]&[1, 4]&98.81&42.83\\
    &50&[0, 30]&[1, 3]&99.45&40.57\\
    &100&[0, 32]&[1, 3]&99.71&44.41\\
             \bottomrule
        \end{tabular} \caption{Performance results: Unit Commitment.}\label{table: performance results UC}
    \end{table}
    
We can see \textcolor{black}{from Figure~\ref{subfig: number of constraints UC}} that the reduced problems generated from the three algorithms have up to $85\%$ fewer constraints than the original MILPs. However, such a reduction results in modest online computational savings when the $\mathcal{B}$-learner + CG is used to train the knn classification model. In contrast, using the $\mathcal{S}$-learner + CG manages to substantially improve the performance of the knn in terms of the online solution time. For instance, \textcolor{black}{we can observe in Table~\ref{table: performance results UC} that }the $\mathcal{B}$-learner + CG with $k=5$ employs around $75\%$ of the time required to solve the original MILP formulation, whereas just $45\%$ of the time is used if the knn is run on the $\mathcal{S}$-learner + CG. This is of particular relevance in real-life applications, such as the UC problem, where the optimality guarantee in a short computational time is a priority.

Indeed, due to the good initialization of the constraint generation method, based on the invariant constraint set, our approach is able to get (almost always) a perfect prediction of the constraints needed to recover the optimal solution of the original MILP formulation\textcolor{black}{, as the third column of Table~\ref{table: performance results UC} shows}. Consequently, the number of iterations of the constraint generation method is $1$ in most of the test instances, and the solution times are reduced\textcolor{black}{, as Figures~\ref{subfig: number of CG iterations UC} and \ref{subfig: online time UC} present}. For instance, \textcolor{black}{Table~\ref{table: performance results UC} states that} our method found the optimal solution in approximately $93\%$ of the instances if $k=5$ is chosen. Moreover, if we observe the results for $k = 100$, nearly $100\%$ of the problems attain the optimal MILP solution in one iteration of the constraint generation method.

To summarize, solving to optimality real-world optimization problems, such as the UC problem, may become a challenge. The presented results have shown that the warm-started constraint generation procedure proposed in this paper, which is based on the invariant constraint set, significantly reduces the computational burden and gets reduced MILPs that are equivalent to the original ones in terms of the optimal solution. Consequently, our approach strengthens and supports the use of machine-learning-aided optimization tools for online applications.

\section{Conclusions and Future Work} \label{sec: Conclusions}
  
Solving MILPs to optimality for online applications using traditional algorithms is not always possible due to their high computational burden. While the literature includes speed-up methods to solve MILP based on machine learning techniques, the obtained solutions may be suboptimal (or even infeasible in certain cases) for the original formulation. In this paper, we propose a machine-learning-aided warm-start constraint generation algorithm that attains the optimal solution of an MILP in a shorter computational time. The proposed approach is based on the offline identification of the invariant constraint sets of previous instances of the target MILP. In doing so, we significantly improve the prediction of invariant constraint sets for unseen instances. Thus, a much smaller number of iterations are needed to run the constraint generation algorithm, and the online computational burden is significantly diminished.

We compare our approach in a synthetic MILP and the unit commitment problem with the traditional constraint generation method and with a warm-started methodology that ignores the information given by the critical non-binding constraints. In both examples, the online computational time is significantly reduced with respect to the comparative strategies. For instance, in our experiments with the unit commitment application, the optimal solution is attained using around $40\%$ of the time needed in the original MILP. This shows the advantage of using our methodology for solving MILP in online applications.
  
In our study, the MILP instances in the training set only differ on the right-hand side of their constraints. While our methodology can also be used with MILPs that are parameterized by the coefficient matrix, further investigation is required to evaluate its performance for this case. Another promising research line is the design of methods for integrating information or expert knowledge on the MILP to be solved into the learning process in order to increase the predictive power of the machine learning engine on the MILP solution. Finally, our approach implies the training of an independent machine learning model for each constraint. In order to improve the performance metrics, a future research line would be to take advantage of the possible relationships among the constraints to train a unique machine learning algorithm.
  
\section*{Acknowledgments}
This work was supported in part by the Spanish Ministry of Science and Innovation through project PID2020-115460GB-I00, in part by the European Research Council (ERC) under the EU Horizon 2020 research and innovation program (grant agreement No. 755705) in part, by the Junta de Andalucía (JA), the Universidad de Málaga (UMA), and the European Regional Development Fund (FEDER) through the research projects P20\_00153 and UMA2018‐FEDERJA‐001, and in part by the Research Program for Young Talented Researchers of the University of M\'alaga under Project B1-2020-15. The authors thankfully acknowledge the computer resources, technical expertise, and assistance provided by the SCBI (Supercomputing and Bioinformatics) center of the University of Málaga.

\newpage
%\bibColoredItems{blue}{lin2017multi, zamfirache2022policy}
\bibliography{\myreferences}

\begin{thebibliography}{22}
\expandafter\ifx\csname natexlab\endcsname\relax\def\natexlab#1{#1}\fi
\providecommand{\url}[1]{\texttt{#1}}
\providecommand{\href}[2]{#2}
\providecommand{\path}[1]{#1}
\providecommand{\DOIprefix}{doi:}
\providecommand{\ArXivprefix}{arXiv:}
\providecommand{\URLprefix}{URL: }
\providecommand{\Pubmedprefix}{pmid:}
\providecommand{\doi}[1]{\href{http://dx.doi.org/#1}{\path{#1}}}
\providecommand{\Pubmed}[1]{\href{pmid:#1}{\path{#1}}}
\providecommand{\bibinfo}[2]{#2}
\ifx\xfnm\relax \def\xfnm[#1]{\unskip,\space#1}\fi
%Type = Article
\bibitem[{Bengio et~al.(2021)Bengio, Lodi \& Prouvost}]{bengio2021machine}
\bibinfo{author}{Bengio, Y.}, \bibinfo{author}{Lodi, A.}, \&
  \bibinfo{author}{Prouvost, A.} (\bibinfo{year}{2021}).
\newblock \bibinfo{title}{Machine learning for combinatorial optimization: A
  methodological tour d’horizon}.
\newblock {\it \bibinfo{journal}{European Journal of Operational Research}\/},
  {\it \bibinfo{volume}{290}\/}, \bibinfo{pages}{405--421}.
  \DOIprefix\doi{10.1016/j.ejor.2020.07.063}.
%Type = Article
\bibitem[{Bertsimas et~al.(2021)Bertsimas, Cory-Wright \&
  Pauphilet}]{bertsimas2021unified}
\bibinfo{author}{Bertsimas, D.}, \bibinfo{author}{Cory-Wright, R.}, \&
  \bibinfo{author}{Pauphilet, J.} (\bibinfo{year}{2021}).
\newblock \bibinfo{title}{A unified approach to mixed-integer optimization
  problems with logical constraints}.
\newblock {\it \bibinfo{journal}{SIAM Journal on Optimization}\/},  {\it
  \bibinfo{volume}{31}\/}, \bibinfo{pages}{2340--2367}.
  \DOIprefix\doi{10.1137/20M1346778}.
%Type = Article
\bibitem[{Bertsimas \& Stellato(2021)}]{bertsimas2021voice}
\bibinfo{author}{Bertsimas, D.}, \& \bibinfo{author}{Stellato, B.}
  (\bibinfo{year}{2021}).
\newblock \bibinfo{title}{The voice of optimization}.
\newblock {\it \bibinfo{journal}{Machine Learning}\/},  {\it
  \bibinfo{volume}{110}\/}, \bibinfo{pages}{249--277}.
  \DOIprefix\doi{10.1007/s10994-020-05893-5}.
%Type = Book
\bibitem[{Bertsimas \& Tsitsiklis(1997)}]{bertsimas1997introduction}
\bibinfo{author}{Bertsimas, D.}, \& \bibinfo{author}{Tsitsiklis, J.}
  (\bibinfo{year}{1997}).
\newblock {\it \bibinfo{title}{Introduction to Linear Optimization}\/}.
\newblock (\bibinfo{edition}{1st} ed.).
\newblock \bibinfo{publisher}{Athena Scientific}.
%Type = Inproceedings
\bibitem[{Bouffard et~al.(2005)Bouffard, Galiana \&
  Arroyo}]{bouffard2005umbrella}
\bibinfo{author}{Bouffard, F.}, \bibinfo{author}{Galiana, F.~D.}, \&
  \bibinfo{author}{Arroyo, J.~M.} (\bibinfo{year}{2005}).
\newblock \bibinfo{title}{Umbrella contingencies in security-constrained
  optimal power flow}.
\newblock In {\it \bibinfo{booktitle}{15th Power systems computation
  conference, PSCC}\/}.
\newblock Volume~\bibinfo{volume}{5}.
%Type = Article
\bibitem[{Calafiore(2010)}]{calafiore2010random}
\bibinfo{author}{Calafiore, G.~C.} (\bibinfo{year}{2010}).
\newblock \bibinfo{title}{Random convex programs}.
\newblock {\it \bibinfo{journal}{SIAM Journal on Optimization}\/},  {\it
  \bibinfo{volume}{20}\/}, \bibinfo{pages}{3427--3464}.
  \DOIprefix\doi{10.1137/090773490}.
%Type = Inproceedings
\bibitem[{Chmiela et~al.(2021)Chmiela, Khalil, Gleixner, Lodi \&
  Pokutta}]{chmiela2021learning}
\bibinfo{author}{Chmiela, A.}, \bibinfo{author}{Khalil, E.~B.},
  \bibinfo{author}{Gleixner, A.}, \bibinfo{author}{Lodi, A.}, \&
  \bibinfo{author}{Pokutta, S.} (\bibinfo{year}{2021}).
\newblock \bibinfo{title}{Learning to schedule heuristics in branch and bound}.
\newblock In \bibinfo{editor}{A.~Beygelzimer}, \bibinfo{editor}{Y.~Dauphin},
  \bibinfo{editor}{P.~Liang}, \& \bibinfo{editor}{J.~W. Vaughan} (Eds.), {\it
  \bibinfo{booktitle}{Advances in Neural Information Processing Systems}\/}.
\newblock \URLprefix \url{https://openreview.net/forum?id=mvEhkIqn45_}.
%Type = Book
\bibitem[{Conforti et~al.(2014)Conforti, Cornuejols \&
  Zambelli}]{conforti2014integer}
\bibinfo{author}{Conforti, M.}, \bibinfo{author}{Cornuejols, G.}, \&
  \bibinfo{author}{Zambelli, G.} (\bibinfo{year}{2014}).
\newblock {\it \bibinfo{title}{Integer Programming}\/}.
\newblock \bibinfo{publisher}{Springer Publishing Company, Incorporated}.
%Type = Book
\bibitem[{Friedman et~al.(2001)Friedman, Hastie \&
  Tibshirani}]{friedman2001elements}
\bibinfo{author}{Friedman, J.}, \bibinfo{author}{Hastie, T.}, \&
  \bibinfo{author}{Tibshirani, R.} (\bibinfo{year}{2001}).
\newblock {\it \bibinfo{title}{The elements of statistical learning}\/}
  Volume~\bibinfo{volume}{1} of {\it \bibinfo{series}{Springer {S}eries in
  {S}tatistics}\/}.
\newblock \bibinfo{publisher}{Springer, Berlin}.
%Type = Article
\bibitem[{Gambella et~al.(2021)Gambella, Ghaddar \&
  Naoum-Sawaya}]{gambella2021optimization}
\bibinfo{author}{Gambella, C.}, \bibinfo{author}{Ghaddar, B.}, \&
  \bibinfo{author}{Naoum-Sawaya, J.} (\bibinfo{year}{2021}).
\newblock \bibinfo{title}{Optimization problems for machine learning: A
  survey}.
\newblock {\it \bibinfo{journal}{European Journal of Operational Research}\/},
  {\it \bibinfo{volume}{290}\/}, \bibinfo{pages}{807--828}.
  \DOIprefix\doi{10.1016/j.ejor.2020.08.045}.
%Type = Book
\bibitem[{G\'omez-Exposito et~al.(2018)G\'omez-Exposito, Conejo \&
  Ca\~nizares}]{gomez2018electric}
\bibinfo{author}{G\'omez-Exposito, A.}, \bibinfo{author}{Conejo, A.~J.}, \&
  \bibinfo{author}{Ca\~nizares, C.} (\bibinfo{year}{2018}).
\newblock {\it \bibinfo{title}{Electric energy systems: analysis and
  operation}\/}.
\newblock \bibinfo{publisher}{CRC press}.
%Type = Book
\bibitem[{Hastie et~al.(2009)Hastie, Tibshirani \&
  Friedman}]{hastie2009elements}
\bibinfo{author}{Hastie, T.}, \bibinfo{author}{Tibshirani, R.}, \&
  \bibinfo{author}{Friedman, J.} (\bibinfo{year}{2009}).
\newblock {\it \bibinfo{title}{The elements of statistical learning: data
  mining, inference, and prediction}\/}.
\newblock \bibinfo{publisher}{Springer Science \& Business Media}.
%Type = Article
\bibitem[{Lin et~al.(2017)Lin, Wang, Zhang, Xiahou \& McDonald}]{lin2017multi}
\bibinfo{author}{Lin, F.}, \bibinfo{author}{Wang, J.}, \bibinfo{author}{Zhang,
  N.}, \bibinfo{author}{Xiahou, J.}, \& \bibinfo{author}{McDonald, N.}
  (\bibinfo{year}{2017}).
\newblock \bibinfo{title}{Multi-kernel learning for multivariate performance
  measures optimization}.
\newblock {\it \bibinfo{journal}{Neural Computing and Applications}\/},  {\it
  \bibinfo{volume}{28}\/}, \bibinfo{pages}{2075--2087}.
  \DOIprefix\doi{10.1007/s00521-015-2164-9}.
%Type = Article
\bibitem[{Lodi et~al.(2020)Lodi, Mossina \& Rachelson}]{lodi2020learning}
\bibinfo{author}{Lodi, A.}, \bibinfo{author}{Mossina, L.}, \&
  \bibinfo{author}{Rachelson, E.} (\bibinfo{year}{2020}).
\newblock \bibinfo{title}{Learning to handle parameter perturbations in
  combinatorial optimization: An application to facility location}.
\newblock {\it \bibinfo{journal}{EURO Journal on Transportation and
  Logistics}\/},  {\it \bibinfo{volume}{9}\/}, \bibinfo{pages}{100023}.
  \DOIprefix\doi{10.1016/j.ejtl.2020.100023}.
%Type = Article
\bibitem[{Minoux(1989)}]{minoux1989networks}
\bibinfo{author}{Minoux, M.} (\bibinfo{year}{1989}).
\newblock \bibinfo{title}{Networks synthesis and optimum network design
  problems: Models, solution methods and applications}.
\newblock {\it \bibinfo{journal}{Networks}\/},  {\it \bibinfo{volume}{19}\/},
  \bibinfo{pages}{313--360}. \DOIprefix\doi{10.1002/net.3230190305}.
%Type = Misc
\bibitem[{OASYS(2022)}]{OASYS2022warm}
\bibinfo{author}{OASYS} (\bibinfo{year}{2022}).
\newblock \bibinfo{title}{{Warm{\_}starting{\_}CG{\_}for{\_}MIO{\_}ML}}.
\newblock
  \bibinfo{howpublished}{\url{https://github.com/groupoasys/Warm_starting_CG_for_MIO_ML}}.
%Type = Article
\bibitem[{{Pineda} et~al.(2020){Pineda}, {Morales} \&
  {Jim\'enez-Cordero}}]{pineda2020data}
\bibinfo{author}{{Pineda}, S.}, \bibinfo{author}{{Morales}, J.~M.}, \&
  \bibinfo{author}{{Jim\'enez-Cordero}, A.} (\bibinfo{year}{2020}).
\newblock \bibinfo{title}{Data-driven screening of network constraints for unit
  commitment}.
\newblock {\it \bibinfo{journal}{IEEE Transactions on Power Systems}\/},  {\it
  \bibinfo{volume}{35}\/}, \bibinfo{pages}{3695--3705}.
  \DOIprefix\doi{10.1109/TPWRS.2020.2980212}.
%Type = Inproceedings
\bibitem[{Taunk et~al.(2019)Taunk, De, Verma \& Swetapadma}]{taunk2019brief}
\bibinfo{author}{Taunk, K.}, \bibinfo{author}{De, S.}, \bibinfo{author}{Verma,
  S.}, \& \bibinfo{author}{Swetapadma, A.} (\bibinfo{year}{2019}).
\newblock \bibinfo{title}{A brief review of nearest neighbor algorithm for
  learning and classification}.
\newblock In {\it \bibinfo{booktitle}{2019 International Conference on
  Intelligent Computing and Control Systems (ICCS)}\/} (pp.
  \bibinfo{pages}{1255--1260}).
\newblock \DOIprefix\doi{10.1109/ICCS45141.2019.9065747}.
%Type = Inbook
\bibitem[{Wolsey(2008)}]{wolsey2008mixed}
\bibinfo{author}{Wolsey, L.~A.} (\bibinfo{year}{2008}).
\newblock \bibinfo{title}{Mixed integer programming}.
\newblock In {\it \bibinfo{booktitle}{Wiley Encyclopedia of Computer Science
  and Engineering}\/} (pp. \bibinfo{pages}{1--10}).
\newblock \bibinfo{publisher}{American Cancer Society}.
\newblock \DOIprefix\doi{10.1002/9780470050118.ecse244}.
%Type = Article
\bibitem[{Xavier et~al.(2021)Xavier, Qiu \& Ahmed}]{xavier2020learning}
\bibinfo{author}{Xavier, A.~S.}, \bibinfo{author}{Qiu, F.}, \&
  \bibinfo{author}{Ahmed, S.} (\bibinfo{year}{2021}).
\newblock \bibinfo{title}{Learning to solve large-scale security-constrained
  unit commitment problems}.
\newblock {\it \bibinfo{journal}{INFORMS Journal on Computing}\/},  {\it
  \bibinfo{volume}{33}\/}, \bibinfo{pages}{739--756}.
  \DOIprefix\doi{10.1287/ijoc.2020.0976}.
%Type = Article
\bibitem[{Yang \& Shami(2020)}]{yang2020on}
\bibinfo{author}{Yang, L.}, \& \bibinfo{author}{Shami, A.}
  (\bibinfo{year}{2020}).
\newblock \bibinfo{title}{On hyperparameter optimization of machine learning
  algorithms: Theory and practice}.
\newblock {\it \bibinfo{journal}{Neurocomputing}\/},  {\it
  \bibinfo{volume}{415}\/}, \bibinfo{pages}{295--316}.
  \DOIprefix\doi{10.1016/j.neucom.2020.07.061}.
%Type = Article
\bibitem[{Zamfirache et~al.(2022)Zamfirache, Precup, Roman \&
  Petriu}]{zamfirache2022policy}
\bibinfo{author}{Zamfirache, I.~A.}, \bibinfo{author}{Precup, R.-E.},
  \bibinfo{author}{Roman, R.-C.}, \& \bibinfo{author}{Petriu, E.~M.}
  (\bibinfo{year}{2022}).
\newblock \bibinfo{title}{Policy iteration reinforcement learning-based control
  using a {G}rey {W}olf optimizer algorithm}.
\newblock {\it \bibinfo{journal}{Information Sciences}\/},  {\it
  \bibinfo{volume}{585}\/}, \bibinfo{pages}{162--175}.
  \DOIprefix\doi{https://doi.org/10.1016/j.ins.2021.11.051}.

\end{thebibliography}

\end{document}